\documentclass[12pt]{article}
\usepackage{multicol}
\usepackage{amsmath, amssymb, amscd, enumerate, amsfonts, latexsym, amsxtra, 
amsthm, bm, multirow}
\usepackage{indentfirst}
\newtheorem{Theorem}{Theorem}[section]
\newtheorem{Proposition}[Theorem]{Proposition}
\newtheorem{Lemma}[Theorem]{Lemma}

\theoremstyle{definition} 
\newtheorem{Remark}[Theorem]{Remark}
\theoremstyle{remark} 
\newtheorem*{Proof}{Proof}

\begin{document}
\allowdisplaybreaks[3]
\centerline{\Large Some results  on certain  finite-dimensional  subalgebras of} \vspace{3mm}
\centerline{\Large the hyperalgebra of   a universal Chevalley group} \vspace{7mm}
\centerline{Yutaka Yoshii 
\footnote{ E-mail address: yutaka.yoshii.6174@vc.ibaraki.ac.jp}}  \vspace{5mm}
\centerline{College of Education,   
Ibaraki University,}
\centerline{2-1-1 Bunkyo, Mito, Ibaraki, 310-8512, Japan}
\begin{abstract}
In the   hyperalgebra of the 
$r$-th Frobenius kernel of  a   
universal Chevalley group over a field of  characteristic $p>0$, we study 
some subsets and the subalgebras generated by them  
and give some results.  We are particularly 
interested in the case that $p$ is 'very small' and the Dynkin diagram is not simply-laced.  
\end{abstract}

\section{Introduction}
Let $\mathbb{F}_p$ be a  prime field of $p$ elements. 
Let $\mathfrak{g}_{\mathbb{C}}$ be a simple complex Lie algebra with 
root vectors $e_{\alpha}$ and $G$ the universal Chevalley group 
of type $\mathfrak{g}_{\mathbb{C}}$ over the algebraic closure $\overline{\mathbb{F}}_p$. Let $\mathcal{U}_{\mathbb{Z}}$ 
be Kostant's $\mathbb{Z}$-form of the universal enveloping algebra 
of $\mathfrak{g}_{\mathbb{C}}$. Then 
the associative $\mathbb{F}_p$-algebra 
$\mathcal{U}= \mathcal{U}_{\mathbb{Z}} \otimes_{\mathbb{Z}} \mathbb{F}_p$ is 
one of the hyperalgebras relative to  $G$. For $r \in \mathbb{Z}_{>0}$, let 
$\mathcal{U}_r$  be the $\mathbb{F}_p$-subalgebra 
of $\mathcal{U}$ generated by all the divided powers 
$e_{\alpha}^{(n)}$ for roots $\alpha$ and  integers $n$ with $0 \leq n \leq p^r-1$. It 
is known that the algebra is related to the $r$-th Frobenius kernel  of $G$.  Let 
$\mathcal{U}^+$  be the $\mathbb{F}_p$-subalgebra of $\mathcal{U}$ 
generated by all $e_{\alpha}^{(n)}$ for positive roots $\alpha$ and nonnegative integers $n$. 
Let $\mathcal{U}_r^+$ be an intersection of $\mathcal{U}^+$ and 
$\mathcal{U}_r$.

Let $\mathcal{V}_r^+$   be the 
$\mathbb{F}_p$-subalgebra of $\mathcal{U}$ generated by 
$1$ and all $e_{\alpha}^{(p^s)}$ for simple roots $\alpha$ and integers $s$ with 
$0 \leq s \leq r-1$. Let $\mathcal{V}_r$ be the 
$\mathbb{F}_p$-subalgebra of $\mathcal{U}$ generated by 
$1$ and all $e_{\alpha}^{(p^s)}$ and $e_{-\alpha}^{(p^s)}$ for simple roots 
$\alpha$ and integers $s$ with $0 \leq s \leq r-1$. 
Then  it is interesting to  study the $\mathbb{F}_p$-subalgebras  
$\mathcal{V}_r^+$ and $\mathcal{V}_r$,  connections between these subalgebras and the 
$\mathbb{F}_p$-subalgebras  $\mathcal{U}_r^+$ and $\mathcal{U}_r$, and 
 generating sets of $\mathcal{U}_r^+$ and $\mathcal{U}_r$. Our goal is to give a series of 
facts about them.    
For example, it turns out that  in most cases 
$\mathcal{V}_r^+$ and $\mathcal{V}_r$ coincide with 
$\mathcal{U}_r^+$ and $\mathcal{U}_r$ respectively, but 
they do not if the type of $\mathfrak{g}_{\mathbb{C}}$ is not simply-laced 
and $p$ is 'very small' (see Theorems \ref{mainthm1} and \ref{mainthm2}). 
On the other hand, it turns out that for each $p$ and each 
type of $\mathfrak{g}_{\mathbb{C}}$, 
the $\mathbb{F}_p$-algebra $\mathcal{U}^+$ is generated by $1$ and all 
$e_{\alpha}^{(p^s)}$ for simple roots $\alpha$ and nonnegative integers $s$, and  
the $\mathbb{F}_p$-algebra $\mathcal{U}$ is generated by $1$ and all  
$e_{\alpha}^{(p^s)}$,  $e_{-\alpha}^{(p^s)}$ for simple roots 
$\alpha$ and nonnegative integers $s$ (see Remark \ref{remark4.10}). 
 
The $\mathbb{F}_p$-algebra $\mathcal{U}$ and its subalgebras are important tools 
on modular representation theory of algebraic groups (for example, see 
\cite{jantzen79} and \cite{jantzen80}). However,  
there seems to be little literature written on the topics in the last paragraph. 
In the quantum case, the  corresponding properties for small quantum groups 
are described in \cite[\S 5]{lusztig90-1} (but when the type is simply-laced). 
But in our case, unlike the quantum one, the algebras $\mathcal{U}_r$ and  
$\mathcal{U}_r^+$ are defined for infinitely many 
$r$'s for each $p$, and so it is important to consider the above topics without excluding 
small $p$'s.

As  preparation for main results, 
in Section 3, we establish  a 
commutation formula of  the elements  $e_{\alpha}^{(n)}$ in $\mathcal{U}_{\mathbb{Z}}^+$ 
for positive roots $\alpha$ and nonnegative integers $n$. 
Then we give main results in Section 4. 
We determine $\mathbb{F}_p$-bases of $\mathcal{V}_r^+$ and a  
minimal generating set  of $\mathcal{U}_r^+$ in Theorem \ref{mainthm1}. As a sequence,  
we also determine $\mathbb{F}_p$-bases of $\mathcal{V}_r$  
and a  generating set  of $\mathcal{U}_r$ in Theorem \ref{mainthm2}.

\section{Preliminaries}
Let $\mathfrak{g}_{\mathbb{C}}$ be a simple complex Lie algebra. Fix a Chevalley basis 
\[\{ e_{\alpha}, h_{\beta} \ |\ \alpha \in \Phi, \beta \in \Delta \}\] 
of $\mathfrak{g}_{\mathbb{C}}$, where $\Phi$ is the set of all roots and 
$\Delta = \{ \alpha_1, \dots, \alpha_{l}\}$ is the set of all simple roots. 
In this paper, we  number the simple roots  following \cite[11.4]{humphreysbook1}.  
Let $\Phi^+$ and $\Phi^-$ be the sets of all positive roots and all negative roots respectively. 
For a simple root $\alpha_i$, we denote $e_{\alpha_i}$, $e_{-\alpha_i}$, and $h_{\alpha_i}$ by 
$e_i$, $f_i$, and $h_i$ respectively. Let $[\cdot, \cdot]$ be the Lie bracket in 
$\mathfrak{g}_{\mathbb{C}}$. 
Then we have $h_i = [e_i,f_i]$ for $i \in \{1, \dots, l \}$. 
For $n \geq 3$ and $z_1, \dots, z_n \in \mathfrak{g}_{\mathbb{C}}$, set 
\[ [z_1, \dots, z_n] = [z_1, [\dots,[z_{n-2}, [z_{n-1},z_n] ] \cdots]].\]

Let $\mathcal{U}_{\mathbb{C}}$ be the universal enveloping algebra of 
$\mathfrak{g}_{\mathbb{C}}$. For $\alpha \in \Phi$ and $n \in \mathbb{Z}_{\geq 0}$, set 
$e_{\alpha}^{(n)} = e_{\alpha}^{n} / n!$. Let $\mathcal{U}_{\mathbb{Z}}$ be the subring of 
$\mathcal{U}_{\mathbb{C}}$ generated by all $e_{\alpha}^{(n)}$ with $\alpha \in \Phi$ 
and $n \in \mathbb{Z}_{\geq 0}$, which is called Kostant's $\mathbb{Z}$-form. Then for 
$i \in \{1, \dots, l\}$, $c \in \mathbb{Z}$, and $n \in \mathbb{Z}_{\geq 0}$, the element 
\[{h_i +c \choose n}= \dfrac{\prod_{k=1}^n (h_i+c-k+1)}{n!} \]
lies in $\mathcal{U}_\mathbb{Z}$. Let $\mathcal{U}_{\mathbb{Z}}^+$, 
$\mathcal{U}_{\mathbb{Z}}^-$, and $\mathcal{U}_{\mathbb{Z}}^0$ be the subrings of 
$\mathcal{U}_{\mathbb{Z}}$ generated by $\left\{ \left. e_{\alpha}^{(n)}\ \right| \ 
\alpha \in \Phi^+, n \in \mathbb{Z}_{\geq 0} \right\}$, 
$\left\{ \left. e_{\alpha}^{(n)}\ \right| \ 
\alpha \in \Phi^-, n \in \mathbb{Z}_{\geq 0} \right\}$, and 
$\left\{ \left. {h_i \choose n}\ \right| \ 
i \in \{ 1, \dots, l \}, n \in \mathbb{Z}_{\geq 0} \right\}$ respectively. Then 
$\mathcal{U}_{\mathbb{Z}}^+$ (resp. $\mathcal{U}_{\mathbb{Z}}^-$) has 
$\left\{ \left. \prod_{\alpha \in \Phi^+} 
 e_{\alpha}^{(n_{\alpha})}\ \right| \ n_{\alpha} \in \mathbb{Z}_{\geq 0} \right\}$ 
(resp. $\left\{ \left. \prod_{\alpha \in \Phi^-} 
 e_{\alpha}^{(n_{\alpha})}\ \right| \ n_{\alpha} \in \mathbb{Z}_{\geq 0} \right\}$) as a 
$\mathbb{Z}$-basis, where an ordering  of 
$\Phi^+$ (resp. $\Phi^-$) in the product is fixed arbitrarily. On the other hand, 
$\mathcal{U}_{\mathbb{Z}}^0$ is commutative and has 
$\left\{ \left. \prod_{i=1}^l {h_i  \choose n_i}\ \right| \ n_i \in \mathbb{Z}_{\geq 0} \right\}$ 
as a $\mathbb{Z}$-basis. The multiplication map 
$\mathcal{U}_{\mathbb{Z}}^- \otimes_{\mathbb{Z}} \mathcal{U}_{\mathbb{Z}}^0 
\otimes_{\mathbb{Z}} \mathcal{U}_{\mathbb{Z}}^+ \rightarrow \mathcal{U}_{\mathbb{Z}}$ 
is a $\mathbb{Z}$-linear isomorphism. The following formulas in 
$\mathcal{U}_{\mathbb{Z}}$ are well-known: 
\[ e_i^{(m)} f_j^{(n)} = f_j^{(n)} e_i^{(m)}, \]
\[ \displaystyle{e_i^{(m)} f_i^{(n)} = \sum_{c=0}^{{\rm min}(m,n)} f_i^{(n-c)} 
{h_i-m-n+2c \choose c} e_i^{(m-c)}}, \] 
\[ \displaystyle{{h_i +s \choose m} e_{\alpha}^{(n)} 
= e_{\alpha}^{(n)} {h_i+s + \langle \alpha, \alpha_i^{\vee} \rangle n \choose m}}, \]
\[ \displaystyle{e_{\alpha}^{(m)} e_{\alpha}^{(n)} = {m+n \choose n} e_{\alpha}^{(m+n)}}, \]
\[ \displaystyle{{h_i \choose m} {h_i \choose n}= \sum_{c=0}^{{\rm min}(m,n) }
{n+m-c \choose m} {m \choose c} {h_i \choose n+m-c}}, \]
\[ \displaystyle{{h_i+s+t \choose m} = 
\sum_{c=0}^{m} {s \choose m-c} {h_i+t \choose c}} \]
for $i,j \in \{1, \dots, l\}\ (i \neq j)$, $m,n \in \mathbb{Z}_{\geq 0}$,   
$s,t \in \mathbb{Z}$, and $\alpha \in \Phi$.  

Let $\mathbb{F}_p$ be a  prime field of $p$ elements. 
Let us define an $\mathbb{F}_p$-algebra $\mathcal{U}$ as the tensor product 
$\mathcal{U}_{\mathbb{Z}} \otimes_{\mathbb{Z}} \mathbb{F}_p$. 
We shall use the same symbols for  images in $\mathcal{U}$ of the elements 
of $\mathcal{U}_{\mathbb{Z}}$ (for example, $e_{\alpha}^{(n)}$, ${h_i +c \choose n}$, 
$[z_1, \dots, z_n]$, and so on).  
Let $\overline{\mathbb{F}}_p$ be an algebraic closure of $\mathbb{F}_p$ and 
$G$  the universal 
Chevalley group of type $\mathfrak{g}_{\mathbb{C}}$ over $\overline{\mathbb{F}}_p$. If 
we regard $G$ as a group scheme defined over $\mathbb{F}_p$,  
$\mathcal{U}$ can be identified with the hyperalgebra of the base change $G_{\mathbb{F}_p}$ 
of $G$ to $\mathbb{F}_p$ (see \cite[I.7.9]{jantzenbook}).   
Let $\mathcal{U}^+$, $\mathcal{U}^-$, and  $\mathcal{U}^0$  
be images in $\mathcal{U}$ of $\mathcal{U}_{\mathbb{Z}}^+$, $\mathcal{U}_{\mathbb{Z}}^-$, 
and $\mathcal{U}_{\mathbb{Z}}^0$ respectively. Then  images in $\mathcal{U}$ of the above 
$\mathbb{Z}$-bases of $\mathcal{U}_{\mathbb{Z}}^+$, 
$\mathcal{U}_{\mathbb{Z}}^-$, and $\mathcal{U}_{\mathbb{Z}}^0$ form  
$\mathbb{F}_p$-bases of $\mathcal{U}^+$, $\mathcal{U}^-$, and $\mathcal{U}^0$ 
respectively. The multiplication map 
$\mathcal{U}^- \otimes_{\mathbb{F}_p} \mathcal{U}^0 
\otimes_{\mathbb{F}_p} \mathcal{U}^+ \rightarrow \mathcal{U}$ 
is an $\mathbb{F}_p$-linear isomorphism. 

For $r \in \mathbb{Z}_{\geq 0}$, let $\mathcal{U}_r$ be the finite-dimensional  
$\mathbb{F}_p$-subalgebra of $\mathcal{U}$ generated by 
$\left\{ \left. e_{\alpha}^{(n)}\ \right| \ 
\alpha \in \Phi, 0\leq n \leq p^r-1 \right\}$. 
If $r>0$, it is the hyperalgebra of the 
$r$-th Frobenius kernel $(G_{\mathbb{F}_p})_r$ of $G_{\mathbb{F}_p}$. 
On the other hand, if $r=0$, we have $\mathcal{U}_0=\mathbb{F}_p$. 
Set 
$\mathcal{U}_r^{+}= \mathcal{U}^{+} \cap \mathcal{U}_r$, 
$\mathcal{U}_r^{-}= \mathcal{U}^{-} \cap \mathcal{U}_r$, and 
$\mathcal{U}_r^{0}= \mathcal{U}^{0} \cap \mathcal{U}_r$. Then 
$\mathcal{U}_r^+$ (resp. $\mathcal{U}_r^-$) has 
$\left\{ \left. \prod_{\alpha \in \Phi^+} 
 e_{\alpha}^{(n_{\alpha})}\ \right| \ 0 \leq n_{\alpha} \leq p^r-1 \right\}$ 
(resp. $\left\{ \left. \prod_{\alpha \in \Phi^-} 
 e_{\alpha}^{(n_{\alpha})}\ \right| \ 0 \leq n_{\alpha} \leq p^r-1 \right\}$) as an 
$\mathbb{F}_p$-basis, where an ordering of $\Phi^+$ (resp. $\Phi^-$) in  
the product is fixed arbitrarily. On the other hand, 
$\mathcal{U}_r^0$ is commutative and has 
$\left\{ \left. \prod_{i=1}^l {h_i  \choose n_i}\ \right| \ 0 \leq n_i \leq p^r-1 \right\}$ 
as an $\mathbb{F}_p$-basis. The multiplication map 
$\mathcal{U}_r^- \otimes_{\mathbb{F}_p} \mathcal{U}_r^0 
\otimes_{\mathbb{F}_p} \mathcal{U}_r^+ \rightarrow \mathcal{U}_r$ 
is an $\mathbb{F}_p$-linear isomorphism. For details, see \cite[II.3]{jantzenbook}. 

There exists an $\mathbb{F}_p$-algebra endomorphism 
${\rm Fr}: \mathcal{U} \rightarrow \mathcal{U}$ defined by 
\[{\rm Fr}\left(e_{\alpha}^{(n)}\right)= 
\left\{ \begin{array}{ll}
{e_{\alpha}^{(n/p)}} & \mbox{if $p\ |\ n$,} \\
0 & \mbox{otherwise} 
\end{array}
\right.,\ \ \ 
{\rm Fr}\left({h_i \choose n}\right)= 
\left\{ \begin{array}{ll}
{h_i \choose n/p} & \mbox{if $p\ |\ n$,} \\
0 & \mbox{otherwise} 
\end{array}
\right.\]
for $a \in \Phi$, $n \in \mathbb{Z}_{\geq 0}$, and $i \in \{ 1, \dots, l\}$.

\section{A commutation formula}

In  this section we give a commutation formula in $\mathcal{U}_{\mathbb{Z}}^+$. Set  
$\nu = |\Phi^+|$ (the number of positive roots). 
Let $W$ be the Weyl group generated by the simple reflections 
$s_i (=s_{\alpha_i})$ with $i \in \{ 1, \dots, l\}$. 
Let $w_0$ be the longest element of $W$ and $s_{i_1} \cdots s_{i_{\nu}}$ its fixed 
reduced expression. Set 
$\beta_1 = \alpha_{i_1}$ and $\beta_k=s_{i_1} \cdots s_{i_{k-1}}(\alpha_{i_k})$ 
for $2 \leq k \leq \nu$. Then we have 
$\Phi^+ = \{ \beta_1, \dots, \beta_{\nu}\}$ 
(see Problem (b) in \cite[Appendix on Finite Reflection Groups II.25]{steinbergbook}). 
We define a total order $\preceq$ in $\Phi^+$ as 
\[ \beta_j \preceq \beta_k \Longleftrightarrow j \leq k. \]
\ 

\begin{Lemma}\label{commform1} 
Let $j,k \in \mathbb{Z}$ 
with $1 \leq j < k \leq \nu$. Then the following hold. \\

\noindent {\rm (i)} For $b \in \mathbb{Z}_{>0}$, the element 
$e_{\beta_k} e_{\beta_j}^{(b)}-e_{\beta_j}^{(b)}e_{\beta_k} $ is $0$ if $k-j=1$ and is 
a $\mathbb{Z}$-linear combination of elements of the form 
$e_{\beta_j}^{(c)} e_{\beta_s}$ with $0 \leq c < b$ and $j < s < k$ if $k-j \geq 2$. \\ 

\noindent {\rm (ii)} For $a \in \mathbb{Z}_{>0}$,  the element 
$e_{\beta_k}^{(a)} e_{\beta_j}-e_{\beta_j}e_{\beta_k}^{(a)} $ is $0$ if $k-j=1$ and is 
a $\mathbb{Z}$-linear combination of elements of the form 
$e_{\beta_s} e_{\beta_k}^{(c)}$ with $0 \leq c < a$ and $j < s < k$ if $k-j \geq 2$. \\ 

\noindent {\rm (iii)} For $a_j, \dots, a_{k-1}  \in \mathbb{Z}_{\geq 0}$, the element 
$e_{\beta_k} e_{\beta_j}^{(a_j)} \cdots e_{\beta_{k-1}}^{(a_{k-1})} -
e_{\beta_j}^{(a_j)} \cdots e_{\beta_{k-1}}^{(a_{k-1})} e_{\beta_k} $ is a $\mathbb{Z}$-linear 
combination of elements of the form $e_{\beta_j}^{(b_j)} \cdots e_{\beta_{k-1}}^{(b_{k-1})}$ 
satisfying the following: \\ \\
$\bullet$ $b_j \leq a_j$. \\
$\bullet$ $\sum_{i=j}^{k-1}b_i \leq \sum_{i=j}^{k-1} a_i$ and if $k-j \geq 2$, then 
$\sum_{i=j+1}^{k-1}b_i \leq \sum_{i=j+1}^{k-1} a_i +1$. \\

\noindent {\rm (iv)} For $a_{j+1}, \dots, a_{k}  \in \mathbb{Z}_{\geq 0}$, the element 
$ e_{\beta_{j+1}}^{(a_{j+1})} \cdots e_{\beta_{k}}^{(a_{k})} e_{\beta_j}-
e_{\beta_j} e_{\beta_{j+1}}^{(a_{j+1})} \cdots e_{\beta_{k}}^{(a_{k})}  $ is a $\mathbb{Z}$-linear 
combination of elements of the form $e_{\beta_{j+1}}^{(b_{j+1})} \cdots e_{\beta_{k}}^{(b_{k})}$ 
satisfying the following: \\ \\
$\bullet$ $b_k \leq a_k$. \\
$\bullet$ $\sum_{i={j+1}}^{k}b_i \leq \sum_{i={j+1}}^{k} a_i$ and if $k-j \geq 2$, then 
$\sum_{i=j+1}^{k-1}b_i \leq \sum_{i=j+1}^{k-1} a_i +1$. \\

\end{Lemma}

\begin{Proof} 
Note that any element in $\mathcal{U}_{\mathbb{Z}}$ is a $\mathbb{Z}$-linear combination 
of elements of the form $\prod_{i=1}^{\nu} e_{\beta_i}^{(n_i)}$ with $n_i \in \mathbb{Z}_{\geq 0}$. 
Since it is ensured that all the coefficients appearing in (i)-(iv) 
are integers, we only have to argue  
the basis elements. We  prove only (i) and (iii) because (ii) and (iv) will follow 
by symmetry. 

We shall show (i). Note that if $\beta_j + \beta_k \in \Phi$, we have 
  $\beta_j \prec \beta_j + \beta_k \prec \beta_k$ (see \cite{papi94}). 
Note also that if $\beta_j + \beta_k \not\in \Phi$, then  
$e_{\beta_k} e_{\beta_j}= e_{\beta_j} e_{\beta_k}$ and hence  
 we have  
$e_{\beta_k} e_{\beta_j}^{(b)}- e_{\beta_j}^{(b)} e_{\beta_k}=0$ in $\mathcal{U}_{\mathbb{Z}}$. 
Since $\beta_j + \beta_{j+1} \not\in \Phi$,  we have nothing to do if $k-j=1$. So we may 
assume that $k-j \geq 2$ and $\beta_j + \beta_k \in \Phi$. 
Then there is an integer $t$ with $j < t < k$ such that 
$\beta_j + \beta_k = \beta_t$. Use induction on 
$b$. If $b=1$,  then we see that 
$e_{\beta_k} e_{\beta_j} - e_{\beta_j}e_{\beta_k}$ is a multiple of 
$e_{\beta_j + \beta_k}= e_{\beta_t}$ by an element in $\mathbb{Z}$. So suppose that 
$b \geq 2$.  In $\mathcal{U}_{\mathbb{C}}$, we have 
\begin{align*}
e_{\beta_k} e_{\beta_j}^{(b)} - e_{\beta_j}^{(b)} e_{\beta_k} 
= \dfrac{1}{b} \left( e_{\beta_k} e_{\beta_j}^{(b-1)} -
e_{\beta_j}^{(b-1)}  e_{\beta_k} \right) e_{\beta_j} + \dfrac{1}{b} e_{\beta_j}^{(b-1)}
\left( e_{\beta_k} e_{\beta_j} - e_{\beta_j} e_{\beta_k} \right).
\end{align*}
The element $e_{\beta_j}^{(b-1)}\left( e_{\beta_k} e_{\beta_j} - e_{\beta_j} e_{\beta_k} \right)$ is 
a multiple of $e_{\beta_j}^{(b-1)} e_{\beta_t}$. 
On the other hand, by induction, the element 
$\left( e_{\beta_k} e_{\beta_j}^{(b-1)} - e_{\beta_j}^{(b-1)} e_{\beta_k} \right)e_{\beta_j}$ is 
a linear combination of  elements of the form 
$e_{\beta_j}^{(c')} e_{\beta_{t'}} e_{\beta_j}$ with $0 \leq c' < b-1$ and $j < t' < k$. Then 
$e_{\beta_j}^{(c')} e_{\beta_{t'}} e_{\beta_j}$ is equal to 
$(c'+1) e_{\beta_j}^{(c'+1)} e_{\beta_{t'}}$ if $\beta_j + \beta_{t'} \not\in \Phi$ and 
to $(c'+1) e_{\beta_j}^{(c'+1)} e_{\beta_{t'}}+c'' e_{\beta_j}^{(c')} e_{\beta_{t''}}$ 
for some $c'' \in \mathbb{Z}$ and $j < t'' < t'$ if 
$\beta_j + \beta_{t'} \in \Phi$. Therefore, (i) is proved. 

We turn to (iii). Use induction on $k-j$. 
It is clear that 
$e_{\beta_k}  e_{\beta_{k-1}}^{(a_{k-1})} - 
 e_{\beta_{k-1}}^{(a_{k-1})} e_{\beta_k}=0$ 
if $k-j=1$. So we may assume that $k-j \geq 2$. 
Moreover, we  may also assume that $\beta_j + \beta_k \in \Phi$ and $a_j >0$. We have 
\begin{align*}
\lefteqn{e_{\beta_k} e_{\beta_j}^{(a_j)} \cdots e_{\beta_{k-1}}^{(a_{k-1})} - 
e_{\beta_j}^{(a_j)} \cdots e_{\beta_{k-1}}^{(a_{k-1})} e_{\beta_k}} \nonumber \\
&= \left( e_{\beta_k} e_{\beta_j}^{(a_j)} - e_{\beta_j}^{(a_j)} e_{\beta_k}\right) 
e_{\beta_{j+1}}^{(a_{j+1})} \cdots e_{\beta_{k-1}}^{(a_{k-1})} + 
e_{\beta_j}^{(a_j)} \left( e_{\beta_k} 
e_{\beta_{j+1}}^{(a_{j+1})} \cdots e_{\beta_{k-1}}^{(a_{k-1})}
-e_{\beta_{j+1}}^{(a_{j+1})} \cdots e_{\beta_{k-1}}^{(a_{k-1})} e_{\beta_k} \right). 
\end{align*}
By induction, the element 
$e_{\beta_k} 
e_{\beta_{j+1}}^{(a_{j+1})} \cdots e_{\beta_{k-1}}^{(a_{k-1})}
-e_{\beta_{j+1}}^{(a_{j+1})} \cdots e_{\beta_{k-1}}^{(a_{k-1})} e_{\beta_k} $ is a linear combination 
of elements of the form $e_{\beta_{j+1}}^{(c_{j+1})} \cdots e_{\beta_{k-1}}^{(c_{k-1})}$ 
satisfying the following: \\ \\
$\bullet$ $c_{j+1} \leq a_{j+1}$. \\ 
$\bullet$ $\sum_{i=j+1}^{k-1}c_{i} \leq \sum_{i=j+1}^{k-1}a_{i}$ and if $k-(j+1) \geq 2$, 
then $\sum_{i=j+2}^{k-1}c_{i} \leq \sum_{i=j+2}^{k-1}a_{i} +1$. \\ \\
Therefore, 
setting $(c_j', c_{j+1}', \dots, c_{k-1}') = (a_j, c_{j+1}, \dots, c_{k-1})$ we see that 
the element 
$$e_{\beta_j}^{(a_j)} \left( e_{\beta_k} 
e_{\beta_{j+1}}^{(a_{j+1})} \cdots e_{\beta_{k-1}}^{(a_{k-1})}
-e_{\beta_{j+1}}^{(a_{j+1})} \cdots e_{\beta_{k-1}}^{(a_{k-1})} e_{\beta_k} \right)$$
is a linear combination 
of elements of the form $e_{\beta_{j}}^{(c_{j}')} \cdots e_{\beta_{k-1}}^{(c_{k-1}')}$ 
satisfying the following: \\ \\
$\bullet$ $c_{j}' = a_{j}$. \\
$\bullet$ $\sum_{i=j}^{k-1}c_{i}' \leq \sum_{i=j}^{k-1}a_{i}$ and if $k-j \geq 2$, 
then $\sum_{i=j+1}^{k-1}c_{i}' \leq \sum_{i=j+1}^{k-1}a_{i} +1$. \\ \\
On the other hand, by (i), 
$e_{\beta_k} e_{\beta_j}^{(a_j)}-e_{\beta_j}^{(a_j)} e_{\beta_k}$ is a linear combination of 
elements of the 
form $e_{\beta_j}^{(c)} e_{\beta_s}$ with $0 \leq c < a_j$ and $j < s < k$. If $s-j \geq 2$, then 
by induction we see that 
\begin{align*}
\lefteqn{e_{\beta_s} e_{\beta_{j+1}}^{(a_{j+1})} \cdots e_{\beta_{s-1}}^{(a_{s-1})} } \nonumber \\
&= \sum_{(d_{j+1}, \dots, d_{s-1}) \in \mathbb{Z}_{\geq 0}^{s-j-1}} 
\eta_1(d_{j+1}, \dots, d_{s-1})
e_{\beta_{j+1}}^{(d_{j+1})} \cdots e_{\beta_{s-1}}^{(d_{s-1})}  +
e_{\beta_{j+1}}^{(a_{j+1})} \cdots e_{\beta_{s-1}}^{(a_{s-1})}  e_{\beta_s},
\end{align*}
where $\eta_1(d_{j+1}, \dots, d_{s-1}) \in \mathbb{Z}$ and 
$(d_{j+1}, \dots, d_{s-1})$ with $\eta_1(d_{j+1}, \dots, d_{s-1}) \neq 0$ satisfies \\ \\
$\bullet$ $d_{j+1} \leq a_{j+1}$. \\
$\bullet$ $\sum_{i=j+1}^{s-1} d_i \leq \sum_{i=j+1}^{s-1} a_i$. \\ \\
Then we have 
\begin{align*}
\lefteqn{e_{\beta_j}^{(c)}e_{\beta_s} e_{\beta_{j+1}}^{(a_{j+1})} \cdots e_{\beta_{k-1}}^{(a_{k-1})}}
\nonumber \\
&= \sum_{(d_{j+1}, \dots, d_{s-1}) \in \mathbb{Z}_{\geq 0}^{s-j-1}} 
\eta_1(d_{j+1}, \dots, d_{s-1})
e_{\beta_j}^{(c)} e_{\beta_{j+1}}^{(d_{j+1})} \cdots e_{\beta_{s-1}}^{(d_{s-1})} 
e_{\beta_s}^{(a_s)} \cdots e_{\beta_{k-1}}^{(a_{k-1})} \nonumber \\
& \ \ + (a_s+1) e_{\beta_j}^{(c)} e_{\beta_{j+1}}^{(a_{j+1})} \cdots e_{\beta_{s-1}}^{(a_{s-1})} 
e_{\beta_s}^{(a_s+1)} e_{\beta_{s+1}}^{(a_{s+1})} \cdots e_{\beta_{k-1}}^{(a_{k-1})}. 
\end{align*}
On the other hand, if $s-j=1$, then 
\begin{align*}
e_{\beta_j}^{(c)}e_{\beta_s} e_{\beta_{j+1}}^{(a_{j+1})} \cdots e_{\beta_{k-1}}^{(a_{k-1})}
= (a_{j+1}+1)e_{\beta_j}^{(c)}e_{\beta_{j+1}}^{(a_{j+1}+1)} 
 e_{\beta_{j+2}}^{(a_{j+2})} \cdots e_{\beta_{k-1}}^{(a_{k-1})}.
\end{align*}
Taking $(d_j', \dots, d_{k-1}')$ as 
\[ \mbox{$(c, d_{j+1}, \dots, d_{s-1}, a_s, \dots, a_{k-1})$\ \ \   or \ \ 
$(c, a_{j+1}, \dots, a_{s-1}, a_s+1, a_{s+1}, \dots, a_{k-1})$} \]
if $s-j \geq 2$ and $(c, a_{j+1}+1, a_{j+2}, \dots, a_{k-1})$ if $s-j=1$, 
we see that the element 
$\left( e_{\beta_k} e_{\beta_j}^{(a_j)} - e_{\beta_j}^{(a_j)} e_{\beta_k}\right) 
e_{\beta_{j+1}}^{(a_{j+1})} \cdots e_{\beta_{k-1}}^{(a_{k-1})}$ is a linear combination of elements 
of the form $e_{\beta_j}^{(d_j')} \cdots e_{\beta_{k-1}}^{(d_{k-1}')}$ satisfying \\ \\
$\bullet$ $d_j' < a_j$. \\
$\bullet$ $\sum_{i=j}^{k-1} d_i' < \sum_{i=j}^{k-1} a_i$ and if $k-j \geq 2$, then 
$\sum_{i=j+1}^{k-1} d_i' \leq  \sum_{i=j+1}^{k-1} a_i +1$. \\ \\ 
Therefore,  (iii) follows. \qed
\end{Proof}
\ 

The commutation formula described in the proposition below will be used in the next 
section. It is a bit more refined than the specialization $v \rightarrow 1$ of  
\cite[Theorem 2.3]{xi99}.   \\

\begin{Proposition}\label{commform2}
For $a,b \in \mathbb{Z}_{>0}$ and  
$j,k \in \mathbb{Z}$ with $1 \leq j < k \leq \nu$, the element 
$e_{\beta_k}^{(a)} e_{\beta_j}^{(b)}-e_{\beta_j}^{(b)}e_{\beta_k}^{(a)}$ 
in $\mathcal{U}_{\mathbb{Z}}$ is a $\mathbb{Z}$-linear 
combination of elements of the form $e_{\beta_j}^{(a_j)} \cdots e_{\beta_k}^{(a_k)}$ 
satisfying the following: \\ \\
$\bullet$ $a_j <b$ and $a_k < a$. \\
$\bullet$ $\sum_{i=j}^{k-1} a_i \leq b$ and $\sum_{i=j+1}^{k}a_i \leq a$. \\

\end{Proposition}

\begin{Proof} 
As in the proof of Lemma \ref{commform1}, it is ensured that 
all the coefficients appearing in the proposition are integers. 
So we only have to show the result 
on basis elements. 
If $a=1$, the proposition follows from Lemma \ref{commform1} (i). Suppose that $a>1$. Then 
\begin{align*}
e_{\beta_k}^{(a)} e_{\beta_j}^{(b)} - e_{\beta_j}^{(b)} e_{\beta_k}^{(a)}  
= \dfrac{1}{a} e_{\beta_k} \left( e_{\beta_k}^{(a-1)} e_{\beta_j}^{(b)} - 
e_{\beta_j}^{(b)} e_{\beta_k}^{(a-1)}\right)  +
\dfrac{1}{a}  \left( e_{\beta_k} e_{\beta_j}^{(b)}  - 
e_{\beta_j}^{(b)} e_{\beta_k}  \right) e_{\beta_k}^{(a-1)}
\end{align*}
(in $\mathcal{U}_{\mathbb{C}}$). By Lemma \ref{commform1} (i), the element 
 $\left( e_{\beta_k} e_{\beta_j}^{(b)} - 
e_{\beta_j}^{(b)} e_{\beta_k} \right) e_{\beta_k}^{(a-1)} $ is 
a linear combination of elements of the 
form $e_{\beta_j}^{(c)}e_{\beta_s} e_{\beta_k}^{(a-1)} $ satisfying 
$0 \leq c <b$  and  $j < s < k$.  
On the other hand, by induction, the element 
$ e_{\beta_k}^{(a-1)} e_{\beta_j}^{(b)} - e_{\beta_j}^{(b)} e_{\beta_k}^{(a-1)}$ is a 
$\mathbb{Z}$-linear combination of elements of the form 
$e_{\beta_j}^{(b_j)} \cdots e_{\beta_k}^{(b_k)}$ satisfying the following: \\ \\
$\bullet$ $b_j < b$ and $b_k <a-1$. \\
$\bullet$ $\sum_{i=j}^{k-1} b_i \leq b$ and $\sum_{i=j+1}^{k} b_i \leq a-1$. \\ \\
Then by Lemma \ref{commform1} (iii), we have 
\begin{align*}
e_{\beta_k} e_{\beta_{j}}^{(b_{j})} \cdots e_{\beta_{k-1}}^{(b_{k-1})} 
= \sum_{(b_{j}' , \dots, b_{k-1}') \in \mathbb{Z}_{\geq 0}^{k-j}} 
\eta_2(b_{j}' , \dots, b_{k-1}') 
e_{\beta_{j}}^{(b_{j}')} \cdots e_{\beta_{k-1}}^{(b_{k-1}')}  + 
 e_{\beta_{j}}^{(b_{j})} \cdots e_{\beta_{k-1}}^{(b_{k-1})}  e_{\beta_k}
\end{align*}
and hence 
\begin{align*}
\lefteqn{e_{\beta_k} e_{\beta_j}^{(b_{j})} \cdots e_{\beta_{k-1}}^{(b_{k-1})} e_{\beta_k}^{(b_k)} } 
\nonumber \\
&= \sum_{(b_{j}' , \dots, b_{k-1}') \in \mathbb{Z}_{\geq 0}^{k-j}} \eta_2(b_{j}' , \dots, b_{k-1}') 
e_{\beta_j}^{(b_{j}')}  \cdots e_{\beta_{k-1}}^{(b_{k-1}')} e_{\beta_k}^{(b_k)} + 
(b_k+1)e_{\beta_j}^{(b_{j})} \cdots e_{\beta_{k-1}}^{(b_{k-1})} e_{\beta_k}^{(b_k+1)},  
\end{align*}
where $\eta_{2}(b_{j}' , \dots, b_{k-1}') \in \mathbb{Z}$ and 
$(b_{j}' , \dots, b_{k-1}')$ with $\eta_{2}(b_{j}' , \dots, b_{k-1}') \neq 0$ 
satisfies the following: \\ \\
$\bullet$ $b_j' \leq b_j$. \\
$\bullet$ $\sum_{i=j}^{k-1} b_i' \leq \sum_{i=j}^{k-1} b_i $. \\ \\ 
Taking $(a_j, \dots, a_k)$ as $(b_{j}', \dots, b_{k-1}', b_k)$ or 
$(b_j, \dots, b_{k-1}, b_k+1)$ we see that the element 
$e_{\beta_k}  \left( e_{\beta_k}^{(a-1)} e_{\beta_j}^{(b)} - 
e_{\beta_j}^{(b)} e_{\beta_k}^{(a-1)}\right) $ is a linear combination of elements 
of the form $e_{\beta_j}^{(a_j)} \cdots e_{\beta_k}^{(a_k)}$ satisfying the following: \\ \\
$\bullet$ $a_j <b$ and $a_k < a$. \\
$\bullet$ $\sum_{i=j}^{k-1} a_i \leq b$ and $\sum_{i=j+1}^k a_i \leq a$. \\ \\ 
Therefore, the proposition is proved. \qed
\end{Proof}
\ 

\begin{Remark}\label{commform2rem} 
Regarding a term $e_{\beta_j}^{(a_j)} \cdots e_{\beta_k}^{(a_k)}$ 
in the proposition, we see the following: \\

\noindent {\rm (i)} If there is an integer $s$ with $j<s<k$ such that $a_s=b$, then 
$e_{\beta_j}^{(a_j)} \cdots e_{\beta_k}^{(a_k)}=e_{\beta_s}^{(b)}e_{\beta_k}^{(a_k)}$. \\

\noindent {\rm (ii)} If there is an integer $s$ with $j<s<k$ such that $a_s=a$, then 
$e_{\beta_j}^{(a_j)} \cdots e_{\beta_k}^{(a_k)}=e_{\beta_j}^{(a_j)}e_{\beta_s}^{(a)}$. \\

\noindent {\rm (iii)} If $b=a$ and there is an integer $s$ with $j<s<k$ such that $a_s=a$, then 
$e_{\beta_j}^{(a_j)} \cdots e_{\beta_k}^{(a_k)}=e_{\beta_s}^{(a)}$ and 
$\beta_s=\beta_j+\beta_k$. \\
\end{Remark}

\section{Main results}
Let $\Lambda$ be the subset of $\mathcal{U}$ consisting of $1$ and all 
$e_i^{(p^{n})}$, $f_i^{(p^{n})}$ with $i \in \{1,\dots, l\}$ and $n \in \mathbb{Z}_{\geq 0}$.  
Set 
\[\Lambda^+ = \Lambda \cap \mathcal{U}^+ = 
\left\{ \left. 1, e_i^{(p^{n})}\ \right| \ i \in \{1,\dots, l\}, n \in \mathbb{Z}_{\geq 0} \right\}, \]
\[\Lambda^- = \Lambda \cap \mathcal{U}^- = 
\left\{ \left. 1, f_i^{(p^{n})}\ \right| \ i \in \{1,\dots, l\}, n \in \mathbb{Z}_{\geq 0} \right\}, \]
\[\Lambda_r = \Lambda \cap \mathcal{U}_r = 
\left\{ \left. 1, e_i^{(p^{n})},   f_i^{(p^{n})} \ \right| \ i \in \{1,\dots, l\}, 
0 \leq n \leq r-1 \right\}, \]
\[\Lambda_r^+ = \Lambda \cap \mathcal{U}_r^+ = 
\left\{ \left. 1, e_i^{(p^{n})} \ \right| \ i \in \{1,\dots, l\}, 
0 \leq n \leq r-1 \right\}, \]
\[\Lambda_r^- = \Lambda \cap \mathcal{U}_r^- = 
\left\{ \left. 1, f_i^{(p^{n})} \ \right| \ i \in \{1,\dots, l\}, 
0 \leq n \leq r-1 \right\} \]
for $r \in \mathbb{Z}_{>0}$. 
Let $\mathcal{V}$,  $\mathcal{V}^+$,  $\mathcal{V}^-$,  $\mathcal{V}_r$, 
 $\mathcal{V}_r^+$, and  $\mathcal{V}_r^-$ be the $\mathbb{F}_p$-subalgebras generated by 
$\Lambda$, $\Lambda^+$, $\Lambda^-$, $\Lambda_r$, $\Lambda_r^+$, and 
$\Lambda_r^-$ respectively. In this section we 
give some results about these $\mathbb{F}_p$-algebras. As preparation, we make some 
observations on the root system $\Phi$. \\ 

Let $\mathbb{E}$ be an euclidean space with inner product $\langle \cdot, \cdot \rangle$ 
spanned by  $\Phi$. 
For $\alpha \in \Phi$, we define the coroot $\alpha^{\vee}$ of $\alpha$ 
as $\alpha^{\vee}= 2\alpha / \langle \alpha, \alpha\rangle$.

For a root $\alpha= \sum_{i=1}^l c_i \alpha_i$ $(c_i \in \mathbb{Z}_{\geq 0})$, the integer 
${\rm ht}(\alpha)=\sum_{i=1}^l c_i$ is called the height of $\alpha$. For two positive roots 
$\alpha, \beta \in \Phi^+$, if there exists a simple root $\alpha_i$ $(i \in \{1, \dots, l\})$ 
such that $\beta = \alpha + \alpha_i$, we  denote this by 
$\alpha \xrightarrow{\alpha_i} \beta$. 
It forms a Hasse diagram $\mathcal{H}(\Phi^+)$ of $\Phi^+$, which 
 is well-known in each type. 

If $\Phi$ is not simply-laced (i.e. $\Phi$ is 
of type ${\rm B}_l\ (l \geq 2)$, ${\rm C}_l\ (l \geq 2)$, ${\rm F}_4$, or ${\rm G}_2$), 
it has two root lengths. Then let $\mathcal{H}(\Phi_{\rm short}^+)$ 
(resp. $\mathcal{H}(\Phi_{\rm long}^+)$) be the subdiagram of $\mathcal{H}(\Phi^+)$ 
whose vertices correspond to all short (resp. long) positive roots. 
Now we shall observe the subdiagrams in 
each type. \\

\noindent $\bullet$ Suppose that $\Phi$ is of type ${\rm B}_l\ (l \geq 2)$. Then 
$\mathcal{H}(\Phi_{\rm long}^+)$ has two connected components whose   
vertices correspond to 
\[ \left\{ \left. \sum_{k=i}^{j} \alpha_k\ \right|\ 1 \leq i \leq j \leq l-1 \right\} 
\mbox{\ \ \ and\ \ \ } 
 \left\{ \left. \sum_{k=i}^j \alpha_k +2 \sum_{k=j+1}^l \alpha_k 
\ \right|\ 1 \leq i \leq j \leq l-1 \right\} \] 
respectively. On the other hand, $\mathcal{H}(\Phi_{\rm short}^+)$ is a connected 
subdiagram as  
\[ \alpha_l \xrightarrow{\alpha_{l-1}} \alpha_{l-1}+\alpha_{l} \xrightarrow{\alpha_{l-2}} \cdots 
\xrightarrow{\alpha_{1}} \alpha_1 + \cdots + \alpha_{l-1}+\alpha_l. \]
\ \\

\noindent $\bullet$ Suppose that $\Phi$ is of type ${\rm C}_l\ (l \geq 2)$. Then 
$\mathcal{H}(\Phi_{\rm long}^+)$ has $l$ connected components and each of them has  
only one vertex corresponding to  
$2\sum_{k=i}^{l-1}  \alpha_k + \alpha_l\ (1 \leq i \leq l) $. On the other hand, 
$\mathcal{H}(\Phi_{\rm short}^+)$ is a connected subdiagram. \\ 

\noindent $\bullet$ Suppose that $\Phi$ is of type ${\rm F}_4$. Then 
$\mathcal{H}(\Phi_{\rm long}^+)$ has four connected components:
\[ \alpha_1 \xrightarrow{\alpha_2} \alpha_1 + \alpha_2 \xleftarrow{\alpha_1} 
\alpha_2, \]
\[ \alpha_2 + 2\alpha_3 \xrightarrow{\alpha_1} \alpha_1 + \alpha_2 +2 \alpha_3 
 \xrightarrow{\alpha_2} \alpha_1+2\alpha_2+2\alpha_3, \]
\[ \alpha_2 + 2\alpha_3 +2\alpha_4 \xrightarrow{\alpha_1} 
\alpha_1 + \alpha_2 +2 \alpha_3+2\alpha_4 
 \xrightarrow{\alpha_2} \alpha_1+2\alpha_2+2\alpha_3+2\alpha_4, \]  
\[ \alpha_1 + 2\alpha_2+4\alpha_3 +2\alpha_4 \xrightarrow{\alpha_2} 
\alpha_1 + 3\alpha_2 +4 \alpha_3+2\alpha_4 
 \xrightarrow{\alpha_1} 2\alpha_1+3\alpha_2+4\alpha_3+2\alpha_4. \]
On the other hand,  $\mathcal{H}(\Phi_{\rm short}^+)$ is a connected subdiagram. \\

\noindent $\bullet$  Suppose that $\Phi$ is of type ${\rm G}_2$. Then 
$\mathcal{H}(\Phi_{\rm long}^+)$ has two connected components: 
\[ \alpha_2\ \ \ \ {\rm and}\ \ \ \ 3\alpha_1+\alpha_2 \xrightarrow{\alpha_2} 
 3\alpha_1+2\alpha_2. \]
On the other hand, $\mathcal{H}(\Phi_{\rm short}^+)$ is a connected 
subdiagram as 
\[ \alpha_1 \xrightarrow{\alpha_2} \alpha_1+\alpha_2 \xrightarrow{\alpha_1}  
2\alpha_1 +\alpha_2. \] 
\ 

From the above  observation, we see the following facts. \\

\begin{Proposition}\label{subdiagram}
Suppose that $\Phi$ is of type ${\rm B}_l\ (l \geq 2)$, 
${\rm C}_l\ (l \geq 2)$, ${\rm F}_4$, or ${\rm G}_2$. 
Then the following hold. \\

\noindent {\rm (i)} $\mathcal{H}(\Phi_{\rm long}^+)$ has a unique 
connected component which contains vertices corresponding to 
long simple roots. \\

\noindent {\rm (ii)}  $\mathcal{H}(\Phi_{\rm short}^+)$ is  connected and contains 
at least one vertex corresponding to a short simple root. \\

\noindent {\rm (iii)} Any connected component of $\mathcal{H}(\Phi_{\rm long}^+)$ 
which does not contain vertices corresponding to  
long simple roots has a unique minimal vertex.  \\

\noindent {\rm (iv)} All the simple roots assigned to arrows in  
$\mathcal{H}(\Phi_{\rm long}^+)$ are long.
\end{Proposition}
\ 

Thus we write $\mathcal{C}_0$ for the unique  component 
in Proposition \ref{subdiagram} (i). Moreover, for a subdiagram 
$\mathcal{S}$ of $\mathcal{H}(\Phi^+)$, by abuse of notation we denote by 
$\mathcal{S}$ again the set of roots assigned to vertices  in the diagram.  \\

\begin{Lemma}\label{root1}
Suppose that $\Phi$ is of type ${\rm B}_l\ (l \geq 2)$, ${\rm C}_l\ (l \geq 2)$, 
or ${\rm F}_4$ and that $\beta$ is 
a long root. Then if there is a root $\alpha$ such that $\beta-\alpha$ is a short root, then 
$\alpha$ is a short root, $\beta-2\alpha$ is a long root, and neither $\beta-3\alpha$ nor  
$\beta + \alpha$ is a root. 
\end{Lemma}

\begin{Proof}
Set $\gamma= \beta-\alpha$. Since
\[2 \langle \gamma, \gamma \rangle =\langle \beta, \beta \rangle 
= \langle \gamma, \gamma \rangle + 2 \langle \gamma, \alpha \rangle + 
\langle \alpha, \alpha \rangle, \] 
we have 
$\langle \gamma, \gamma \rangle / \langle \alpha, \alpha \rangle 
= \langle \gamma, \alpha^{\vee} \rangle +1$. 
Since $\gamma$ is short and $\langle \gamma, \alpha^{\vee} \rangle$ is an integer, 
we must have $\langle \alpha, \alpha \rangle=\langle \gamma, \gamma \rangle $ and hence 
$\alpha$ is short and $\langle \gamma, \alpha^{\vee} \rangle=0$. Moreover, 
since $\langle \gamma-\alpha, \gamma-\alpha \rangle/
\langle \alpha, \alpha \rangle =2$, $\gamma-\alpha$ must be a long root. Finally, since 
$\langle \beta-3\alpha, \alpha^{\vee} \rangle =-4$ and 
$\langle \beta+\alpha, \alpha^{\vee} \rangle =4$, it follows that neither $\beta-3\alpha$ nor  
$\beta+ \alpha$ is a root. \qed
\end{Proof}
\ 

\begin{Remark}\label{root1rem} 
Under the assumption  of the lemma, the $\alpha$-string through $\beta$ is 
\[ \beta-2\alpha ,\ \beta-\alpha ,\  
\beta. \]
Moreover, by \cite[Theorem 25.2 (d)]{humphreysbook1} 
we see that  $[e_{\alpha}, e_{\beta-\alpha}]= 2 c_1 e_{\beta}$, 
$[e_{\alpha}, e_{\beta-2\alpha}]=  c_2 e_{\beta-\alpha}$, and 
\[ e_{\alpha}^{(2)}e_{\beta-2\alpha}- e_{\beta-2\alpha}e_{\alpha}^{(2)}
=c_1 c_2 e_{\beta} + c_2 e_{\beta-\alpha} e_{\alpha} = 
-c_1 c_2 e_{\beta} + c_2 e_{\alpha} e_{\beta-\alpha}  \] in 
$\mathcal{U}_{\mathbb{Z}}$ for some $c_1,c_2 \in \{ \pm 1\}$. \\
\end{Remark}

\begin{Lemma}\label{root2}
Suppose that $\Phi$ is of type ${\rm B}_l\ (l \geq 2)$, ${\rm C}_l\ (l \geq 2)$, 
or ${\rm F}_4$ and that $\alpha$ is 
a long root satisfying $\alpha= \beta + \gamma$ for some roots $\beta, \gamma \in \Phi$. 
Then $\beta$ and $\gamma$ have the same length. 
\end{Lemma} 

\begin{Proof}
Suppose that $\beta$ is long. Then the equality 
$\langle \alpha, \alpha \rangle = \langle \alpha-\gamma, \alpha-\gamma \rangle $ implies 
$\langle \alpha, \gamma^{\vee} \rangle =1$. On the other hand, the inequality 
$\langle \alpha-\gamma, \alpha-\gamma \rangle \geq 
\langle \gamma, \gamma \rangle $ implies $\langle \gamma, \alpha^{\vee} \rangle \leq 1$. 
Thus we must have $\langle \gamma, \alpha^{\vee} \rangle = 1$. Therefore, we see that 
$\langle \alpha, \alpha \rangle = \langle \gamma, \gamma \rangle $ and hence that 
$\gamma$ is long. 

Suppose that $\beta$ is short. Then the equality 
$\langle \alpha, \alpha \rangle = 2\langle \alpha-\gamma, \alpha-\gamma \rangle $ 
implies $\langle \alpha, \gamma^{\vee} \rangle = 
\langle \alpha, \alpha \rangle/2\langle \gamma, \gamma \rangle +1$. Therefore, we see that 
$\langle \alpha, \alpha \rangle/2\langle \gamma, \gamma \rangle$ is an integer and 
hence that $\gamma$ is short. \qed
\end{Proof}
\ 

\begin{Lemma}\label{root3}
Suppose that $\Phi$ is of type ${\rm B}_l\ (l \geq 2)$, ${\rm C}_l\ (l \geq 2)$, 
or ${\rm F}_4$ and that $\alpha$ is 
a long root satisfying $\alpha= \beta + \gamma$ for some $\beta, \gamma \in \Phi^+$. 
If $\alpha \not\in \mathcal{C}_0$ and 
$\beta, \gamma \in \mathcal{C}_0 \cup \mathcal{H}(\Phi_{\rm short}^+)$, then   
both $\beta$ and $\gamma$ are short. 
\end{Lemma}

\begin{Proof}
Let $\delta \in \Phi^+$. Then 
by Proposition \ref{subdiagram} (iv) we see the following: 
\[ \delta \in \mathcal{C}_0 \Longleftrightarrow 
\mbox{$\delta$ is a sum of some long simple roots}. \]
It follows from Lemma \ref{root2} 
that $\beta$ and $\gamma$ must have the same length. 
Suppose that both $\beta$ and $\gamma$ are long. 
Then since $\beta, \gamma \in \mathcal{C}_0$, they must be sums of some long simple roots.   
Thus we have $\alpha \in \mathcal{C}_0$, which is contradiction. \qed
\end{Proof}
\ 

Now we study a connection between $\mathcal{V}_r^+$ and $\mathcal{U}_r^+$. \\

\begin{Proposition}\label{subalggeneral}
Suppose that $\Phi$ satisfies one of the following: \\ 

\noindent $\bullet$ $p$ is arbitrary and $\Phi$ is of type 
${\rm A}_l\ (l \geq 1)$,  ${\rm D}_l\ (l \geq 4)$, or ${\rm E}_l\ (l  \in \{6,7,8\})$. 

\noindent $\bullet$ $p \geq 3$  and $\Phi$ is of type 
${\rm B}_l\ (l \geq 2)$,  ${\rm C}_l\ (l \geq 2)$, or ${\rm F}_4$. 

\noindent $\bullet$ $p \geq 5$  and $\Phi$ is of type 
${\rm G}_2$. \\ 

\noindent Then $\mathcal{V}_r^+ =\mathcal{U}_r^+ $. 
\end{Proposition} 

\begin{Proof}
Let $\alpha \in \Phi^+$. Consider a path 
\begin{align}\label{path1}
\alpha_{i_1} = \gamma_1 \xrightarrow{\alpha_{i_2}} \gamma_2 
\xrightarrow{\alpha_{i_3}} \cdots \xrightarrow{\alpha_{i_n}} \gamma_n = \alpha
\end{align}
in $\mathcal{H}(\Phi^+)$, where $n= {\rm ht}(\alpha)$ and 
$i_1, \dots, i_n \in \{1, \dots, l\}$. 

Suppose that $r=1$. We would like to show that $e_{\alpha} \in \mathcal{V}_1^+$. 
It is clear  when $n=1$, so  assume that $n>1$. Consider the equality  
$d e_{\alpha}= [e_{i_n}, e_{\gamma_{n-1}}]$ in $\mathcal{U}_{\mathbb{Z}}$,  where  
$d \in \mathbb{Z}$. 
Use \cite[Theorem 25.2 (d)]{humphreysbook1}. 
If $\Phi$ is of type ${\rm A}_l$, ${\rm D}_l$, or ${\rm E}_l$, then since 
$\gamma_{n-1} -\alpha_{i_n} \not\in \Phi$,  $d$ must lie in $\{ \pm 1\}$. 
If $\Phi$ is of type ${\rm B}_l$, ${\rm C}_l$, or ${\rm F}_4$, then since 
$\gamma_{n-1}-a\alpha_{i_n} \in \Phi$ with $a \in \mathbb{Z}_{\geq 0}$ implies 
$a \in \{ 0,1\}$, 
$d$ must lie in $\{ \pm 1, \pm 2\}$. If $\Phi$ is of type ${\rm  G}_2$, then 
since $\gamma_{n-1}-a\alpha_{i_n}\in \Phi$ with 
$a \in \mathbb{Z}_{\geq 0}$ implies $a \in \{ 0,1,2\}$, 
$d$ must lie in $\{ \pm 1, \pm 2, \pm 3\}$. Thus the assumption in the proposition  implies 
that $d \neq 0$ in $\mathbb{F}_p$. In $\mathcal{U}$, since $e_{\gamma_{n-1}} \in \mathcal{V}_1^+$ by 
induction on $n$, we have $e_{\alpha} \in \mathcal{V}_1^+$. Therefore, we obtain 
$\mathcal{V}_1^+=\mathcal{U}_1^+$. 

Suppose that $r>1$. Since $\mathcal{V}_{r-1}^+ =\mathcal{U}_{r-1}^+ $ by induction on $r$, 
it is enough to check that $e_{\alpha}^{(p^{r-1})} \in \mathcal{V}_r^+$. 
Consider  path (\ref{path1}) and use induction on $n$. It is clear for $n=1$, 
so suppose that $n >1$. 
By Proposition \ref{commform2},  we must have 
\[e_{i_n}^{(p^{r-1})} e_{\gamma_{n-1}}^{(p^{r-1})}- 
e_{\gamma_{n-1}}^{(p^{r-1})}e_{i_n}^{(p^{r-1})} = c e_{\alpha}^{(p^{r-1})} + z \]
in $\mathcal{U}$,  where $c \in \mathbb{F}_p$ and $z \in \mathcal{U}_{r-1}^+$.  
Applying the map ${\rm Fr}^{r-1}$ to the equality we have 
$[e_{i_n}, e_{\gamma_{n-1}}] = c e_{\alpha}$, and hence $c$ must be equal to the above $d$ 
in $\mathbb{F}_p$. Since 
$e_{\gamma_{n-1}}^{(p^{r-1})} \in \mathcal{V}_r^+$ by induction and since 
$d \neq 0$ in $\mathbb{F}_p$, we have $e_{\alpha}^{(p^{r-1})} \in \mathcal{V}_r^+$. 
Therefore, we obtain $\mathcal{V}_{r}^+ =\mathcal{U}_{r}^+ $.   \qed
\end{Proof}
\ 

In turn, we have to deal with the cases that $\Phi$ is of type 
${\rm B}_l\ (l \geq 2)$, ${\rm C}_l\ (l \geq 2)$, or ${\rm F}_4$ when $p=2$ and 
$\Phi$ is of type ${\rm G}_2$ when $p=2$ or $3$. 
In type ${\rm G}_2$, we will often carry out some direct calculations. 
So for simplicity, we denote the positive root vectors 
$e_{\alpha_1}$, $e_{\alpha_2}$, $e_{\alpha_1+\alpha_2}$, $e_{2\alpha_1+\alpha_2}$, 
$e_{3\alpha_1+\alpha_2}$, and $e_{3\alpha_1+2\alpha_2}$ by 
$e_1$, $e_2$, $e_{12}$, $e_{112}$, $e_{1112}$, and $e_{11122}$ respectively. 
Without loss of generality, we may assume that these elements are chosen such that 
\begin{align}\label{commformg1}
e_{12}= e_1e_2-e_2e_1, 
\end{align} 
\begin{align}\label{commformg2}
2e_{112}= e_1e_{12}-e_{12}e_1, 
\end{align} 
\begin{align}\label{commformg3}
3e_{1112}= e_1e_{112}-e_{112}e_1, 
\end{align} 
\begin{align}\label{commformg4}
e_{11122}= e_2e_{1112}-e_{1112}e_2, 
\end{align} 
\begin{align}\label{commformg5}
3e_{11122}= e_{112}e_{12}-e_{12}e_{112} 
\end{align} 
in $\mathfrak{g}_{\mathbb{C}}$ (for example, see \cite[\S 5]{lusztig90-2}). 
Then for $a,b \in \mathbb{Z}_{>0}$, we have 
\begin{align}\label{commformg6}
e_1^{(a)}  e_2^{(b)} = \sum_{(t_1,t_2,t_3,t_4,t_5,t_6) \in \mathcal{A}_1} 
e_2^{(t_1)} e_{12}^{(t_2)} e_{11122}^{(t_3)} e_{112}^{(t_4)} e_{1112}^{(t_5)} e_1^{(t_6)},
\end{align} 
\begin{align}\label{commformg7}
e_1^{(a)}  e_{112}^{(b)} = \sum_{(t_1,t_2,t_3) \in \mathcal{A}_2} 
3^{t_2} e_{112}^{(t_1)} e_{1112}^{(t_2)} e_1^{(t_3)},
\end{align} 
\begin{align}\label{commformg8}
e_{1112}^{(a)}  e_2^{(b)} = \sum_{(t_1,t_2,t_3) \in \mathcal{A}_2} 
(-1)^{t_2} e_{2}^{(t_1)} e_{11122}^{(t_2)} e_{1112}^{(t_3)},
\end{align} 
in $\mathcal{U}_{\mathbb{Z}}$, 
where
\[ \mathcal{A}_1 = \{ (t_1, \dots, t_6) \in \mathbb{Z}_{\geq 0}^6\ |\ 
t_1+t_2+2t_3+t_4+t_5=b, t_2+3t_3+2t_4+3t_5+t_6=a\}, \]
\[ \mathcal{A}_2 = \{ (t_1, t_2, t_3) \in \mathbb{Z}_{\geq 0}^3\ |\ 
t_1+t_2=b, t_2+t_3=a\}. \]
\

\begin{Proposition}\label{subalgBCF1}
Suppose that $p=2$ and $\Phi$ is of type 
${\rm B}_l\ (l \geq 2)$, ${\rm C}_l\ (l \geq 2)$, or ${\rm F}_4$. Then the following hold. \\

\noindent {\rm (i)} $\mathcal{V}_{r}^+$ contains $\mathcal{U}_{r-1}^+$. \\

\noindent {\rm (ii)} Let $\alpha \in \Phi^+$. Then 
$\alpha \in \mathcal{C}_0 \cup \mathcal{H}(\Phi_{\rm short}^+)$ if and only if 
$e_{\alpha}^{(2^{r-1})} \in \mathcal{V}_r^+$.  
\end{Proposition}

\begin{Proof}
Let $\alpha \in \Phi^+$. Consider a path 
\begin{align}\label{path2}
\alpha_{i_1} = \gamma_1 \xrightarrow{\alpha_{i_2}} \gamma_2 
\xrightarrow{\alpha_{i_3}} \cdots \xrightarrow{\alpha_{i_n}} \gamma_n = \alpha
\end{align}
in $\mathcal{H}(\Phi^+)$, where $n= {\rm ht}(\alpha)$. 

Suppose that $r=1$. Then since $\mathcal{U}_0^+=\mathbb{F}_p$, (i) is clear. 
We shall show (ii). Suppose that  
$\alpha \in \mathcal{C}_0 \cup \mathcal{H}(\Phi_{\rm short}^+)$. Then we can choose  
the roots assigned to vertices in path (\ref{path2}) such that all of them have the same length 
(see Proposition \ref{subdiagram} (i) and (ii)). Use induction on $n$. Now we write 
$d e_{\alpha} = [e_{i_n}, e_{\gamma_{n-1}}]$ for some $d \in \mathbb{Z}$. 
Then by \cite[Theorem 25.2 (d)]{humphreysbook1} 
$d$ must be $\pm 1$. Indeed, the fact that $\gamma_{n-1}$ and 
$\alpha=\gamma_{n-1}+\alpha_{i_n}$ have the same length 
implies $\langle \gamma_{n-1}, \alpha_{i_n}^{\vee} \rangle = -1$, 
and hence $\gamma_{n-1}-\alpha_{i_n}$ is not a root. In $\mathcal{U}$, since 
$e_{\gamma_{n-1}} \in \mathcal{V}_1^+$ by induction, 
we obtain $e_{\alpha} \in \mathcal{V}_1^+$. 
On the other hand,  suppose that  
$\alpha \not\in \mathcal{C}_0 \cup \mathcal{H}(\Phi_{\rm short}^+)$. 
We wish to show that 
$e_{\alpha} \not\in \mathcal{V}_1^+$. Suppose that $e_{\alpha} \in \mathcal{V}_1^+$. 
Then there is a path as  (\ref{path2}) 
such that $e_{\gamma_k}=[e_{i_k}, e_{\gamma_{k-1}}]$ in $\mathcal{U}$ for any integer $k$ 
with $2 \leq k \leq n$ (and $e_{\alpha} = [e_{i_n}, \dots, e_{i_1}]$ in $\mathcal{U}$). 
It follows from the assumption of $\alpha$ that there is an integer $m$ with $2 \leq m \leq n$ 
such that $\gamma_{m-1}$ is short and 
$\gamma_{m}(=\gamma_{m-1}+\alpha_{i_m})$ is long. Then the fact that 
$\langle \gamma_{m}, \gamma_{m} \rangle / \langle \gamma_{m-1}, \gamma_{m-1} \rangle=2$ 
implies that $\langle \gamma_{m-1}, \alpha_{i_m}^{\vee} \rangle =0$ and $\alpha_{i_m}$ is short. 
Thus we see that  
$\langle \gamma_{m-1}-\alpha_{i_m}, \gamma_{m-1}-\alpha_{i_m} \rangle
/\langle \gamma_{m-1}, \gamma_{m-1} \rangle =2$ and hence that 
$\gamma_{m-1}-\alpha_{i_m}$ is a (long positive) root. On the other hand, since 
$\langle \gamma_{m-1}-2\alpha_{i_m}, \alpha_{i_m}^{\vee} \rangle =-4$, 
$\gamma_{m-1}-2\alpha_{i_m}$ is not a root. Thus by 
\cite[Theorem 25.2 (d)]{humphreysbook1} we have 
$[e_{i_m}, e_{\gamma_{m-1}}] =  2 c e_{\gamma_m}$ in $\mathcal{U}_{\mathbb{Z}}$ 
for some $c \in \{ \pm 1\}$. 
But since $p=2$, we have $[e_{i_m}, e_{\gamma_{m-1}}]  = 0$ in $\mathcal{U}$, which is 
contradiction. Therefore, we must have $e_{\alpha} \not\in \mathcal{V}_1^+$ and hence (ii)  
for $r=1$ is proved. 

From now on suppose that $r>1$. We first show  
(i). Since    
$\mathcal{V}_{r-1}^+$ contains $\mathcal{U}_{r-2}^+$ by induction on $r$, it is enough to 
check that $e_{\alpha}^{(2^{r-2})} \in \mathcal{V}_r^+$. We check 
it by considering a path as (\ref{path2}) and using induction on 
$n(={\rm ht}(\alpha))$. It is 
clear for $n=1$ and hence we may assume that $n>1$. 
Suppose first that $\langle \gamma_n, \gamma_n \rangle \leq 
\langle \gamma_{n-1}, \gamma_{n-1} \rangle $. Then since 
$\langle \gamma_{n-1}-\alpha_{i_n}, \alpha_{i_n}^{\vee} \rangle \leq -3$, 
$\gamma_{n-1}-\alpha_{i_n}$ is not a root. It follows 
from \cite[Theorem 25.2 (d)]{humphreysbook1} that there exists $d \in \{ \pm 1\}$ 
such that $d e_{\alpha}=[e_{i_n}, e_{\gamma_{n-1}}]$ in $\mathcal{U}_{\mathbb{Z}}$. 
Thus by Proposition \ref{commform2} we see that 
\[ e_{i_n}^{(2^{r-2})} e_{\gamma_{n-1}}^{(2^{r-2})}- 
e_{\gamma_{n-1}}^{(2^{r-2})}e_{i_n}^{(2^{r-2})}=
d e_{\alpha}^{(2^{r-2})} + z_1\]
in $\mathcal{U}$, where $z_1 \in \mathcal{U}_{r-2}^+$. Since 
$e_{\gamma_{n-1}}^{(2^{r-2})} \in \mathcal{V}_r^+$ by induction on $n$, we obtain 
$e_{\alpha}^{(2^{r-2})} \in \mathcal{V}_r^+$. 
In turn, suppose that $\langle \gamma_n, \gamma_n \rangle > 
\langle \gamma_{n-1}, \gamma_{n-1} \rangle $. Then a similar argument to the previous 
paragraph shows that  
$\gamma_{n-1}-\alpha_{i_n} \in \Phi^+$, so we can choose $i_{n-1}$ such that 
$i_{n}=i_{n-1}$. We may assume that $\alpha_{i_n} \succ \alpha-2\alpha_{i_n}$ because 
a similar argument works even when $\alpha_{i_n} \prec \alpha-2\alpha_{i_n}$. Now 
Proposition \ref{commform2} implies that 
\[e_{i_n}^{(2^{r-1})} e_{\alpha-2\alpha_{i_n}}^{(2^{r-2})}- 
e_{\alpha-2\alpha_{i_n}}^{(2^{r-2})} e_{i_n}^{(2^{r-1})} = 
c_1 e_{\alpha}^{(2^{r-2})} +c_2 e_{\alpha-\alpha_{i_n}}^{(2^{r-2})} e_{i_n}^{(2^{r-2})} +z_2\]
in $\mathcal{U}$ for some $z_2 \in \mathcal{V}_{r-1}^+$ and 
$c_1, c_2 \in \mathbb{F}_2$. 
By applying the map ${\rm Fr}^{r-2}$ to the equality we have 
\[ e_{i_n}^{(2)} e_{\alpha-2\alpha_{i_n}} -e_{\alpha-2\alpha_{i_n}}  e_{i_n}^{(2)} 
= c_1 e_{\alpha} +c_2 e_{\alpha-\alpha_{i_n}} e_{i_n} +{\rm Fr}^{r-2}(z_2)\]
in $\mathcal{U}$. Then  by Remark \ref{root1rem} 
we see that $c_1, c_2 \in \{ \pm 1\}$ (and ${\rm Fr}^{r-2}(z_2)=0$).   
Since $e_{\alpha-2\alpha_{i_n}}^{(2^{r-2})}$ and $e_{\alpha-\alpha_{i_n}}^{(2^{r-2})}$ lie in  
$\mathcal{V}_r^+$ by induction, we obtain $e_{\alpha}^{(2^{r-2})} \in 
\mathcal{V}_r^+$. 
Therefore, we have shown that $\mathcal{V}_r^+$ contains 
$\mathcal{U}_{r-1}^+$ and (i) is proved.    

Let us show (ii) (for $r>1$). 
Suppose that $\alpha \in \mathcal{C}_0 \cup \mathcal{H}(\Phi_{\rm short}^+)$. We  
would like to show that $e_{\alpha}^{(2^{r-1})} \in \mathcal{V}_r^+$. In  path 
(\ref{path2}), we choose all the roots assigned to vertices 
such that they have the same length. Use 
induction on $n(={\rm ht}(\alpha))$. It is clear for $n=1$, so we assume that $n>1$. 
Then there exists $d \in \{ \pm 1 \}$ such that 
$d e_{\alpha}= [e_{i_n}, e_{\gamma_{n-1}}]$ in $\mathcal{U}_{\mathbb{Z}}$ as in 
the second paragraph. Thus we have
\[ e_{i_n}^{(2^{r-1})} e_{\gamma_{n-1}}^{(2^{r-1})} - 
e_{\gamma_{n-1}}^{(2^{r-1})} e_{i_n}^{(2^{r-1})} =
d e_{\alpha}^{(2^{r-1})}+z_3\]
in $\mathcal{U}$ for some $z_3 \in \mathcal{U}_{r-1}^+$ by Proposition \ref{commform2}. 
Since $e_{\gamma_{n-1}}^{(2^{r-1})} \in \mathcal{V}_r^+$ by induction on $n$, we obtain 
$e_{\alpha}^{(2^{r-1})} \in \mathcal{V}_r^+$.  

Finally, suppose that $\alpha \not\in \mathcal{C}_0 \cup \mathcal{H}(\Phi_{\rm short}^+)$. 
If $e_{\alpha}^{(2^{r-1})} \in \mathcal{V}_r^+$, we have 
$e_{\alpha} = {\rm Fr}^{r-1}\left(e_{\alpha}^{(2^{r-1})}\right) \in \mathcal{V}_1^+$, 
which  contradicts  the result for $r=1$. 
Therefore, we obtain $e_{\alpha}^{(2^{r-1})} \not\in \mathcal{V}_r^+$ and the proof is 
complete. \qed
\end{Proof}
\ 

\begin{Proposition}\label{subalgBCF2}
Suppose that $p=2$ and $\Phi$ is of type 
${\rm B}_l\ (l \geq 2)$, ${\rm C}_l\ (l \geq 2)$, or ${\rm F}_4$. 
Let $\mathcal{C}_1, \dots, \mathcal{C}_t$ $(t \geq 1)$ 
be the connected components of $\mathcal{H}(\Phi_{\rm long}^+)$ not containing 
vertices corresponding to long 
simple roots and $\theta_i$  a unique  minimal root in $\mathcal{C}_i$ for $1 \leq i \leq t$. 
Then the subset $\Lambda_r^+ \cup 
\left\{ \left. e_{\theta_i}^{(2^{r-1})}\ \right| \ 1 \leq i \leq t \right\}$ in $\mathcal{U}$ 
generates the $\mathbb{F}_2$-subalgebra $\mathcal{U}_r^+$.
\end{Proposition}

\begin{Proof}
Let $\widehat{\mathcal{V}}_r^+$ be the $\mathbb{F}_2$-subalgebra generated by the subset 
$\Lambda_r^+ \cup \left\{ \left. e_{\theta_i}^{(2^{r-1})}\ \right| \ 1 \leq i \leq t \right\}$ 
in $\mathcal{U}$.  We wish to show that $\widehat{\mathcal{V}}_r^+= \mathcal{U}_r^+$. 
Thanks to Proposition \ref{subalgBCF1}, it is enough to show that for any 
$\alpha \in \mathcal{C}_1 \cup \cdots \cup \mathcal{C}_t$, the element 
$e_{\alpha}^{(2^{r-1})}$ lies in  
$\widehat{\mathcal{V}}_r^+$. 
Assume that $\alpha \in \mathcal{C}_i$. Then there is a path 
\begin{align}\label{path3}
\theta_{i} = \gamma_1 \xrightarrow{\alpha_{i_2}} \gamma_2 
\xrightarrow{\alpha_{i_3}} \cdots \xrightarrow{\alpha_{i_m}} \gamma_m = \alpha
\end{align}
in $\mathcal{C}_i$ for some $m \in \mathbb{Z}_{>0}$, 
 where all $\gamma_j$'s are long. It is clear that 
$\alpha=\theta_i$ and $e_{\alpha} = e_{\theta_i} \in \widehat{\mathcal{V}}_r^+$ 
when $m=1$, so assume that $m>1$. 
Since the root $\gamma_{m-1}= \gamma_m - \alpha_{i_m}$ is long, we see that 
$\gamma_m - 2\alpha_{i_m} (= \gamma_{m-1} - \alpha_{i_m})$ is not a root. Indeed, 
$\langle \gamma_m, \gamma_m \rangle = \langle \gamma_{m-1}, \gamma_{m-1} \rangle$ 
implies that $\langle \gamma_m-2\alpha_{i_m}, \alpha_{i_m}^{\vee} \rangle =-3$. 
Thus by \cite[Theorem 25.2 (d)]{humphreysbook1} 
there is an integer $c \in \{ \pm 1 \}$  such that 
$[e_{i_m}, e_{\gamma_{m-1}}] = c e_{\gamma_m}$ in $\mathcal{U}_{\mathbb{Z}}$. 
Then Proposition \ref{commform2} implies that 
\[ e_{i_m}^{(2^{r-1})} e_{\gamma_{m-1}}^{(2^{r-1})} - 
e_{\gamma_{m-1}}^{(2^{r-1})} e_{i_m}^{(2^{r-1})} = c e_{\gamma_m}^{(2^{r-1})}+z \]
in $\mathcal{U}$, where $z \in \mathcal{U}_{r-1}^+$. Since 
$e_{\gamma_{m-1}}^{(2^{r-1})} \in \widehat{\mathcal{V}}_r^+$ by  induction 
on $m$, we see that 
$e_{\alpha}^{(2^{r-1})} =e_{\gamma_{m}}^{(2^{r-1})}  \in \widehat{\mathcal{V}}_r^+$. 
Therefore, we obtain $\widehat{\mathcal{V}}_r^+=\mathcal{U}_r^+$. \qed
\end{Proof}
\ 

\begin{Proposition}\label{subalgG}
Suppose that $\Phi$ is of type ${\rm G}_2$. Then the following hold. \\

\noindent {\rm (i)} If $p=2$, then $\mathcal{V}_r^+$ contains $\mathcal{U}_{r-1}^+$, 
$e_{12}^{(2^{r-1})}$ lies in $\mathcal{V}_r^+$, and 
$e_{112}^{(2^{r-1})}, e_{1112}^{(2^{r-1})}$, 
and $e_{11122}^{(2^{r-1})}$ do not lie in 
$\mathcal{V}_r^+$. \\ 

\noindent {\rm (ii)} If $p=3$, then $\mathcal{V}_r^+$ contains $\mathcal{U}_{r-1}^+$, 
$e_{12}^{(3^{r-1})}$ and $e_{112}^{(3^{r-1})}$ 
lie  in $\mathcal{V}_r^+$, and 
$e_{1112}^{(3^{r-1})}$ and $e_{11122}^{(3^{r-1})}$ 
do not lie in 
$\mathcal{V}_r^+$. 
\end{Proposition}

\begin{Proof}
Assume that $p=2$ or $3$. 
It is clear from equalities (\ref{commformg1})-(\ref{commformg5}) 
that $e_{12} \in \mathcal{V}_1^+$ and 
$e_{112}, e_{1112}, e_{11122} \not\in \mathcal{V}_1^+$ when $p=2$, and 
$e_{12}, e_{112} \in \mathcal{V}_1^+$ and 
$e_{1112}, e_{11122} \not\in \mathcal{V}_1^+$ when $p=3$. So assume that $r>1$. 
It follows from induction on $r$ that $\mathcal{V}_{r-1}^+ \supseteq \mathcal{U}_{r-2}^+ $. 

Choosing $a=b=p^{r-2}$ in  equality (\ref{commformg6}), we see that 
\[ e_1^{(p^{r-2})} e_2^{(p^{r-2})}= e_2^{(p^{r-2})} e_1^{(p^{r-2})} +e_{12}^{(p^{r-2})} +z_1\]
in $\mathcal{U}$, where $z_1 \in \mathcal{U}_{r-2}^+$. Thus we obtain 
$ e_{12}^{(p^{r-2})} \in \mathcal{V}_r^+$. 

Choosing $a=2p^{r-2}$ and $b=p^{r-2}$ in  equality (\ref{commformg6}), we see that 
\[ e_1^{(2p^{r-2})} e_2^{(p^{r-2})}= e_2^{(p^{r-2})} e_1^{(2p^{r-2})} +e_{112}^{(p^{r-2})} 
+e_{12}^{(p^{r-2})} e_1^{(p^{r-2})}  +z_2\]
in $\mathcal{U}$, where $z_2 \in \mathcal{V}_{r-1}^+$. Thus we obtain 
$ e_{112}^{(p^{r-2})} \in \mathcal{V}_r^+$.

Choosing $a=3p^{r-2}$ and $b=p^{r-2}$ in  equality (\ref{commformg6}), we see that 
\[ e_1^{(3p^{r-2})} e_{2}^{(p^{r-2})}= 
e_{2}^{(p^{r-2})} e_1^{(3p^{r-2})} + e_{12}^{(p^{r-2})} e_1^{(2p^{r-2})} + 
e_{112}^{(p^{r-2})} e_1^{(p^{r-2})} +e_{1112}^{(p^{r-2})} 
+z_3\]
in $\mathcal{U}$, where $z_3 \in \mathcal{V}_{r}^+$. Thus we obtain 
$ e_{1112}^{(p^{r-2})} \in \mathcal{V}_r^+$. 

Choosing $a=b=p^{r-2}$ in  equality (\ref{commformg8}), we see that 
\[ e_{1112}^{(p^{r-2})} e_{2}^{(p^{r-2})}= e_{2}^{(p^{r-2})} e_{1112}^{(p^{r-2})} -e_{11122}^{(p^{r-2})} 
+z_4\]
in $\mathcal{U}$, where $z_4 \in \mathcal{U}_{r-2}^+$. Thus we obtain 
$ e_{11122}^{(p^{r-2})} \in \mathcal{V}_r^+$. 
Therefore, we obtain $\mathcal{V}_{r}^+ \supseteq \mathcal{U}_{r-1}^+ $. 

Choosing $a=b=p^{r-1}$ in  equality (\ref{commformg6}), we see that 
\[ e_1^{(p^{r-1})} e_2^{(p^{r-1})}= e_2^{(p^{r-1})} e_1^{(p^{r-1})} +e_{12}^{(p^{r-1})} +z_5\]
in $\mathcal{U}$, where $z_5 \in \mathcal{U}_{r-1}^+$. Thus we obtain 
$ e_{12}^{(p^{r-1})} \in \mathcal{V}_r^+$. 

Choosing $a=2p^{r-1}$ and $b=p^{r-1}$ in  equality (\ref{commformg6}), we see that 
\[ e_1^{(2p^{r-1})} e_2^{(p^{r-1})}= e_2^{(p^{r-1})} e_1^{(2p^{r-1})} +
e_{112}^{(p^{r-1})} + e_{12}^{(p^{r-1})} e_1^{(p^{r-1})}+z_6\]
in $\mathcal{U}$, where $z_6 \in \mathcal{V}_{r}^+$. Thus if $p=3$, then 
$e_1^{(2 \cdot 3^{r-1})} $ 
lies in $\mathcal{V}_{r}^+$ and hence we obtain $e_{112}^{(3^{r-1})}  \in \mathcal{V}_{r}^+$. 

On the other hand, since  
$e_{\gamma} ={\rm Fr}^{r-1}\left(e_{\gamma}^{(p^{r-1})}\right) 
\not\in \mathcal{V}_1^+$, we must have 
$e_{\gamma}^{(p^{r-1})} \not\in \mathcal{V}_{r}^+$ for 
$\gamma \in \{ 2\alpha_1+\alpha_2, 3\alpha_1+\alpha_2, 3\alpha_1+2\alpha_2\}$ when $p=2$ 
and $\gamma \in \{ 3\alpha_1+\alpha_2, 3\alpha_1+2\alpha_2\}$ when $p=3$. \qed
\end{Proof}
\ 

\begin{Remark}\label{remark4.10} It follows from Propositions 
\ref{subalggeneral}, \ref{subalgBCF1}, and \ref{subalgG} that 
$\mathcal{U}_{n-1}^+ \subseteq \mathcal{V}_n^+$ and similarly 
$\mathcal{U}_{n-1}^- \subseteq \mathcal{V}_n^-$ for any 
$n \in \mathbb{Z}_{> 0}$ in every case. This fact implies 
that $\mathcal{V}^+ =\mathcal{U}^+$, $\mathcal{V}^- =\mathcal{U}^-$, and  
$\mathcal{V} =\mathcal{U}$ in every case. 
\end{Remark}
\ 

The following theorem is one of the main results in this paper. \\

\begin{Theorem}\label{mainthm1}
$\mathcal{V}_r^+$ has 
$\left\{ \left. \prod_{\alpha \in \Phi^+} e_{\alpha}^{(n_{\alpha})}\ \right| \  
0 \leq n_{\alpha} \leq p^{a_{\alpha}}-1\right\}$ as an $\mathbb{F}_p$-basis, where  
an ordering of $\Phi^+$  in the product is fixed arbitrarily, and 
$\Lambda_r^+ \cup \left\{ \left. e_{\alpha}^{(p^{r-1})}\ \right| \ 
\alpha \in \Theta \right\}$ is a minimal generating set of the  
$\mathbb{F}_p$-subalgebra $\mathcal{U}_r^+$ of $\mathcal{U}$. Here $a_{\alpha}$ and $\Theta$ are defined as follows: \\

\noindent {\rm (a)} If one of the following conditions holds, then $a_{\alpha} =r$ for any 
$\alpha \in \Phi^+$
and $\Theta$ is empty: \\ \\
$\bullet$ $p$ is arbitrary and $\Phi$ is of type ${\rm A}_l\ (l \geq 1)$, 
${\rm D}_l\ (l \geq 4)$, or 
${\rm E}_l\ (l \in \{6,7,8\})$. \\
$\bullet$ $p \geq 3$ and $\Phi$ is of type ${\rm B}_l\ (l \geq 2)$, 
${\rm C}_l\ (l \geq 2)$, or ${\rm F}_4$. \\
$\bullet$ $p \geq 5$ and $\Phi$ is of type ${\rm G}_2$. \\

\noindent {\rm (b)} If $p=2$ and $\Phi$ is of type ${\rm B}_l\ (l \geq 2)$, then 
$\Theta=\{\alpha_{l-1}+2\alpha_l\}$ and $a_{\alpha} = r-1$ if 
$\alpha$ is of the form 
$\sum_{k=i}^{j-1} \alpha_k +2 \sum_{k=j}^{l} \alpha_k$ $(1 \leq i < j \leq l)$ 
and $a_{\alpha}=r$ otherwise. \\

\noindent {\rm (c)} If $p=2$ and $\Phi$ is of type ${\rm C}_l\ (l \geq 2)$, then 
$\Theta=\left\{ \left. 2\sum_{k=i}^{l-1} \alpha_k+\alpha_l
\ \right|\ 1 \leq i \leq l-1\right\}$ 
and $a_{\alpha} = r-1$ if 
$\alpha$ is of the form 
$2\sum_{k=i}^{l-1} \alpha_k+\alpha_l$ $(1 \leq i \leq l-1)$ 
and $a_{\alpha}=r$ otherwise. \\

\noindent {\rm (d)}  If $p=2$ and $\Phi$ is of type ${\rm F}_4$, then 
$\Theta=\{\alpha_2+2\alpha_3, \alpha_2+2\alpha_3+2\alpha_4 \}$ 
and 
$a_{\alpha}=r-1$ if $\alpha$ lies in 
\[ \left\{ \begin{array}{l}
\alpha_2+2\alpha_3,\ \  \alpha_1+\alpha_2+2\alpha_3,\ \ \alpha_1+2\alpha_2+2\alpha_3,\ \  
\alpha_2+2\alpha_3+2\alpha_4,  \\
\alpha_1+\alpha_2+2\alpha_3+2\alpha_4,\ \  \alpha_1+2\alpha_2+2\alpha_3+2\alpha_4,\ \  \alpha_1+2\alpha_2+4\alpha_3+2\alpha_4, \\
 \alpha_1+3\alpha_2+4\alpha_3+2\alpha_4,\ \  
2\alpha_1+3\alpha_2+4\alpha_3+2\alpha_4 
\end{array}\right\} \] 
and $a_{\alpha}=r$ otherwise. \\

\noindent {\rm (e)}  If $p=2$ and $\Phi$ is of type ${\rm G}_2$, then 
$\Theta=\{2\alpha_1+\alpha_2\}$ and $a_{\alpha}=r-1$ if $\alpha$ lies in 
$\{2\alpha_1+\alpha_2, 3\alpha_1+\alpha_2, 3\alpha_1+2\alpha_2\}$ and 
$a_{\alpha}=r$ otherwise. \\

\noindent {\rm (f)}  If $p=3$ and $\Phi$ is of type ${\rm G}_2$, then 
$\Theta=\{3\alpha_1+\alpha_2\}$ and $a_{\alpha}=r-1$ if $\alpha$ lies in 
$\{3\alpha_1+\alpha_2, 3\alpha_1+2\alpha_2\}$ and 
$a_{\alpha}=r$ otherwise. 
\end{Theorem}

\begin{Proof}
Note that $\Lambda_r^+$ is a minimal generating set of the $\mathbb{F}_p$-algebra 
$\mathcal{V}_r^+$. 
So all the claims under (a) follow from Proposition \ref{subalggeneral}. 

We first prove that $\Lambda_r^+ \cup \left\{ \left. e_{\alpha}^{(p^{r-1})}\ \right| \ 
\alpha \in \Theta \right\}$ generates $\mathcal{U}_r^+$ under (b)-(f). 
The  claim under (b) or (c) follows from  Proposition \ref{subalgBCF2}. 
We shall show the  claim under (d). 
Suppose that $p=2$ and $\Phi$ is of type ${\rm F}_4$. 
Let $\widehat{\mathcal{V}}_r'^+$ be the 
$\mathbb{F}_2$-subalgebra of $\mathcal{U}$ generated by 
$\Lambda_r^+ \cup \left\{ \left. e_{\alpha}^{(2^{r-1})}\ \right| \ 
\alpha \in \Theta \right\}$. 
In this situation we only have to show that 
$e_{\alpha_1+2\alpha_2+4\alpha_3+2\alpha_4}^{(2^{r-1})} \in \widehat{\mathcal{V}}_r'^+$ and 
then the  claim (i.e. $\widehat{\mathcal{V}}_r'^+=\mathcal{U}_r^+$) 
follows from Proposition \ref{subalgBCF2}. By Proposition 
\ref{commform2} we have  
\[e_1^{(2^{r-1})} e_{\alpha_2+2\alpha_3}^{(2^{r-1})} = 
e_{\alpha_2+2\alpha_3}^{(2^{r-1})} e_1^{(2^{r-1})} + 
c_1e_{\alpha_1+\alpha_2+2\alpha_3}^{(2^{r-1})}+ z_1,\]
\[e_{\alpha_1+\alpha_2+2\alpha_3}^{(2^{r-1})} e_{\alpha_2+2\alpha_3+2\alpha_4}^{(2^{r-1})} = 
e_{\alpha_2+2\alpha_3+2\alpha_4}^{(2^{r-1})} e_{\alpha_1+\alpha_2+2\alpha_3}^{(2^{r-1})} + 
c_2e_{\alpha_1+2\alpha_2+4\alpha_3+2\alpha_4}^{(2^{r-1})}+ z_2\]
in $\mathcal{U}$, where $c_1,c_2 \in \mathbb{F}_2$ and $z_1, z_2 \in \mathcal{U}_{r-1}^+$. 
Applying the map ${\rm Fr}^{r-1}$ to these equalities, we have 
\[e_1 e_{\alpha_2+2\alpha_3} = 
e_{\alpha_2+2\alpha_3} e_1 + 
c_1e_{\alpha_1+\alpha_2+2\alpha_3},\]  
\[e_{\alpha_1+\alpha_2+2\alpha_3} e_{\alpha_2+2\alpha_3+2\alpha_4} = 
e_{\alpha_2+2\alpha_3+2\alpha_4} e_{\alpha_1+\alpha_2+2\alpha_3} + 
c_2e_{\alpha_1+2\alpha_2+4\alpha_3+2\alpha_4}\] 
in $\mathcal{U}$. 
By \cite[Theorem 25.2 (d)]{humphreysbook1}, we must have $c_1=c_2=1$. 
Therefore, we see that  $e_{\alpha_1+\alpha_2+2\alpha_3}^{(2^{r-1})}$ 
and $e_{\alpha_1+2\alpha_2+4\alpha_3+2\alpha_4}^{(2^{r-1})}$ lie in 
$\widehat{\mathcal{V}}_r'^+$ (note that $\widehat{\mathcal{V}}_r'^+$ contains $\mathcal{U}_{r-1}^+$), as required.    
Now we shall show the above claim under (e) or (f).   
Suppose that $p=2$ or $3$ and $\Phi$ is of type ${\rm G}_2$. Let 
$\widehat{\mathcal{V}}_r^+$ be the $\mathbb{F}_p$-subalgebra generated by 
$\Lambda_r^+ \cup \left\{ e_{112}^{(2^{r-1})} \right\}$ (hence by 
$\mathcal{V}_r^+ \cup \left\{ e_{112}^{(2^{r-1})} \right\}$) if $p=2$ and by 
$\Lambda_r^+ \cup \left\{ e_{1112}^{(3^{r-1})} \right\}$ (hence by 
$\mathcal{V}_r^+ \cup \left\{ e_{1112}^{(3^{r-1})} \right\}$) if $p=3$. 
Thanks to Proposition \ref{subalgG},  to show that 
$\widehat{\mathcal{V}}_r^+=\mathcal{U}_r^+$, it suffices to check that each of 
$e_{1112}^{(2^{r-1})}$ and $e_{11122}^{(2^{r-1})}$ for $p=2$ and  
$e_{11122}^{(3^{r-1})}$ for $p=3$ lies in $\widehat{\mathcal{V}}_r^+$. 
Choosing $a=b=p^{r-1}$ in  equality (\ref{commformg7}), we see that 
\[ e_1^{(p^{r-1})} e_{112}^{(p^{r-1})}= 
e_{112}^{(p^{r-1})} e_1^{(p^{r-1})} +3^{p^{r-1}}e_{1112}^{(p^{r-1})} 
+z_3\]
in $\mathcal{U}$, where $z_3 \in \mathcal{U}_{r-1}^+$. Thus we obtain 
$ e_{1112}^{(2^{r-2})} \in \widehat{\mathcal{V}}_r^+$ when $p=2$. 
Choosing $a=b=p^{r-1}$ in  equality (\ref{commformg8}), we see that 
\[ e_{1112}^{(p^{r-1})} e_{2}^{(p^{r-1})}= e_{2}^{(p^{r-1})} e_{1112}^{(p^{r-1})} -e_{11122}^{(p^{r-1})} 
+z_4\]
in $\mathcal{U}$, where $z_4 \in \mathcal{U}_{r-1}^+$. Thus we obtain 
$ e_{11122}^{(p^{r-1})} \in \widehat{\mathcal{V}}_r^+$, as required. 

We next prove that the generating set 
$\Lambda_r^+ \cup \left\{ \left. e_{\alpha}^{(p^{r-1})}\ \right| \ 
\alpha \in \Theta \right\}$ is minimal under (b)-(f). 
Since $\mathcal{V}_r^+ \neq \mathcal{U}_r^+$ in these cases, 
the minimality of 
$\Lambda_r^+ \cup \left\{ \left. e_{\alpha}^{(p^{r-1})}\ \right| \ 
\alpha \in \Theta \right\}$ is clear if 
$\Theta$ is a singleton set. Thus we may assume that $p=2$ and 
$\Phi$ is of type ${\rm C}_l\ (l \geq 3)$ or ${\rm F}_4$. Let $\beta \in \Theta$. 
We only have to show that 
$e_{\beta}^{(2^{r-1})}$ is not contained in the $\mathbb{F}_2$-subalgebra generated by 
$\Lambda_r^+ \cup \left\{ \left. e_{\alpha}^{(2^{r-1})}\ \right| \ 
\alpha \in \Theta, \alpha \neq \beta \right\}$. Thanks to the map 
${\rm Fr}^{r-1}$, it is enough 
to check that $e_{\beta}$ is not contained in the $\mathbb{F}_2$-subalgebra of 
$\mathcal{U}$ generated by $\{ e_{\alpha}\ |\ \alpha \in \Phi^+, \alpha \neq \beta \}$. 
So suppose that $e_{\beta}$ lies in the $\mathbb{F}_2$-subalgebra. 
Then there must exist two positive 
roots $\gamma_1$ and $\gamma_2$ such that $\beta= \gamma_1+ \gamma_2$ and 
$e_{\beta} = [e_{\gamma_1}, e_{\gamma_2}]$ in $\mathcal{U}_1^+$. But in this situation, 
it is easy to see that both $\gamma_1$ and $\gamma_2$ must be short. Thus by  
Remark  \ref{root1rem} we have $[e_{\gamma_1}, e_{\gamma_2}]=0$ in 
$\mathcal{U}_1^+$, which is contradiction. Therefore, the minimality is proved.  

It remains now to show the first claim under (b)-(f). Note that the dimension of 
the $\mathbb{F}_p$-span of the subset $\left\{ \left. \prod_{\alpha \in \Phi^+} e_{\alpha}^{(n_{\alpha})}\ \right| \  
0 \leq n_{\alpha} \leq p^{a_{\alpha}}-1\right\}$ in $\mathcal{U}$ is independent of 
a choice of a fixed total order in $\Phi^+$. So it is enough to show the claim for  a 
fixed order in $\Phi^+$ and hence we write $\Phi^+= \{\beta_1, \dots, \beta_{\nu}\}$ 
($\nu=|\Phi^+|$) as in Section 3. For convenience, 
if $\Phi$ is of type ${\rm G}_2$,  set 
\[\beta_1=\alpha_2, \ \beta_2=\alpha_1+\alpha_2,\ \beta_3=3\alpha_1+2\alpha_2,\ 
\beta_4=2\alpha_1+\alpha_2, \ \beta_5=3\alpha_1+\alpha_2,\ \beta_6=\alpha_1. \]
It gives a total ordering of  $\Phi^+$ as in Section 3 relative to $w_0= s_2s_1s_2s_1s_2s_1$. 

All we need is to show  the following 
proposition. 
\end{Proof} 
\ 

\begin{Proposition}
Suppose that one of the following holds. \\ \\ 
$\bullet$ $p=2$ and $\Phi$ is of type 
${\rm B}_l\ (l \geq 2)$, ${\rm C}_l\ (l \geq 2)$, or ${\rm F}_4$. \\
$\bullet$ $p=2$ or $3$ and $\Phi$ is of type  ${\rm G}_2$. \\

\noindent Moreover, set $a_i=a_{\beta_i}$ for $1 \leq i \leq \nu$ and let $\widetilde{\mathcal{V}}_r^+$ 
be the $\mathbb{F}_p$-span of the subset 
\[\left\{ \left. \prod_{i=1}^{\nu} e_{\beta_i}^{(n_i)}\ \right| \  
0 \leq n_i \leq p^{a_i}-1\right\}\] 
in $\mathcal{U}$. Then 
$\widetilde{\mathcal{V}}_r^+=\mathcal{U}_r^+$.   
\end{Proposition}

\begin{Proof}
Since $\widetilde{\mathcal{V}}_r^+$ contains $\Lambda_r^+$, it is enough to show that 
$\widetilde{\mathcal{V}}_r^+$ is closed under multiplication. For that, it is enough to show the 
following:  
For $1 \leq j < k \leq \nu$ and $0 \leq b_i \leq p^{a_i}-1$ with $i \in \{j,k\}$, if we write 
the element $e_{\beta_k}^{(b_k)} e_{\beta_j}^{(b_j)}$ in $\mathcal{U}$ as an 
$\mathbb{F}_p$-linear combination of elements of the form 
$e_{\beta_j}^{(c_j)} \cdots e_{\beta_k}^{(c_k)}$ with $c_i \in \mathbb{Z}_{\geq 0}$ 
$(j \leq i \leq k)$, then any such term satisfies $0 \leq c_i \leq p^{a_i}-1$ for each 
$i \in \{j, \dots, k\}$. So consider the element $e_{\beta_k}^{(b_k)} e_{\beta_j}^{(b_j)}$ and 
let $\mathcal{V}_{(r,j,k)}$ be the $\mathbb{F}_p$-span of all 
$e_{\beta_j}^{(c_j)} \cdots e_{\beta_k}^{(c_k)}$ with $0 \leq c_i \leq p^{a_i}-1$ 
$(j \leq i \leq k)$. 
We have to show that $e_{\beta_k}^{(b_k)} e_{\beta_j}^{(b_j)} \in 
\mathcal{V}_{(r,j,k)}$. 
As preparation, for $n \in \mathbb{Z}_{\geq 0}$, let $\mathcal{V}_{(n,j,k)}'$ be 
the $\mathbb{F}_p$-span of all $e_{\beta_j}^{(c_j)} \cdots e_{\beta_k}^{(c_k)}$ with $0 \leq c_i \leq p^{n}-1$ $(j \leq i \leq k)$. 
Proposition \ref{commform2} implies that $\mathcal{V}_{(n,j,k)}'$ is an 
$\mathbb{F}_p$-subalgebra of $\mathcal{U}_n^+$. Since each $a_i$ is $r-1$ or $r$,   
$\mathcal{V}_{(r,j,k)}$ contains $\mathcal{V}_{(r-1,j,k)}'$.

Suppose that  $b_k \leq p^{r-1}-1$. Then by Proposition \ref{commform2} we have 
\[ e_{\beta_k}^{(b_k)}  e_{\beta_j}^{(b_j)} =  e_{\beta_j}^{(b_j)}  e_{\beta_k}^{(b_k)} 
+ \sum_{(d_j, \dots, d_k) \in \mathbb{Z}^{k-j+1}} \rho(d_j, \dots, d_k) 
e_{\beta_j}^{(d_j)}  \cdots e_{\beta_k}^{(d_k)}\]
in $\mathcal{U}$, where $\rho(d_j, \dots, d_k) \in \mathbb{F}_p$ and 
$(d_j, \dots, d_k)$ with 
$\rho(d_j, \dots, d_k) \neq 0$ satisfies the following: \\ 
 
\noindent $\bullet$ $d_j < b_j$. \\
$\bullet$ $\sum_{i=j+1}^k d_{i} \leq b_k (\leq p^{r-1}-1)$. \\

\noindent So $e_{\beta_k}^{(b_k)}  e_{\beta_j}^{(b_j)} $ must 
lie in $\mathcal{V}_{(r,j,k)}$. Similarly,   $e_{\beta_k}^{(b_k)}  e_{\beta_j}^{(b_j)} $ 
also lies in $\mathcal{V}_{(r,j,k)}$ even when $b_j \leq p^{r-1}-1$. 

From now on, suppose that $p^{r-1} \leq b_j,b_k \leq p^{r}-1$. Then $a_j=a_k=r$.  Since 
$e_{\beta_k}^{(b_k)}  e_{\beta_j}^{(b_j)} = e_{\beta_j}^{(b_j)} e_{\beta_k}^{(b_k)} \in 
\mathcal{V}_{(r,j,k)}$ if $\beta_j + \beta_k \not\in \Phi$, we may assume that 
$\beta_j + \beta_k \in \Phi$. Then there is an integer $s$ with $j < s < k$ such that 
$\beta_j + \beta_k =\beta_s$. If we write 
$b_i = p^{r-1} +b_i'$ with $0 \leq b_i' \leq p^{r-1}-1$ for $i \in \{ j,k\}$, we have 
$e_{\beta_k}^{(b_k)}  e_{\beta_j}^{(b_j)} =e_{\beta_k}^{(b_k')}  e_{\beta_k}^{(p^{r-1})} 
e_{\beta_j}^{(p^{r-1})}  e_{\beta_j}^{(b_j')}  $ 
in $\mathcal{U}$. By Remark \ref{commform2rem} we see that 
\begin{align}\label{eqinprop1}
e_{\beta_k}^{(p^{r-1})}  e_{\beta_j}^{(p^{r-1})} = e_{\beta_j}^{(p^{r-1})} e_{\beta_k}^{(p^{r-1})} 
+ d e_{\beta_s}^{(p^{r-1})} +z_1 
\end{align}
in $\mathcal{U}$ for some $d \in \mathbb{F}_p$ and $z_1 \in \mathcal{V}_{(r-1,j,k)}'$. 
Then we have 
\begin{align}\label{eqinprop2}
e_{\beta_k}^{(b_k)}  e_{\beta_j}^{(b_j)} 
=  e_{\beta_k}^{(b_k')}  e_{\beta_j}^{(p^{r-1})} e_{\beta_k}^{(p^{r-1})}  e_{\beta_j}^{(b_j')} +
d e_{\beta_k}^{(b_k')}  e_{\beta_s}^{(p^{r-1})}   e_{\beta_j}^{(b_j')}  + e_{\beta_k}^{(b_k')}  z_1  e_{\beta_j}^{(b_j')}.  
\end{align}
Clearly, the summand $e_{\beta_k}^{(b_k')}  z_1  e_{\beta_j}^{(b_j')}$ 
lies in $\mathcal{V}_{(r-1,j,k)}'$. In turn, since 
\[ e_{\beta_k}^{(b_k')}  e_{\beta_j}^{(p^{r-1})} = e_{\beta_j}^{(p^{r-1})} e_{\beta_k}^{(b_k')} +z_2 \]
and 
\[ e_{\beta_k}^{(p^{r-1})}  e_{\beta_j}^{(b_j')} = e_{\beta_j}^{(b_j')} e_{\beta_k}^{(p^{r-1})} +z_3, \] 
where $z_2,z_3 \in \mathcal{V}_{(r-1,j,k)}'$, $e_{\beta_k}^{(b_k')}  e_{\beta_j}^{(p^{r-1})} e_{\beta_k}^{(p^{r-1})}  e_{\beta_j}^{(b_j')} $ is an $\mathbb{F}_p$-linear combination of elements of the form 
$ e_{\beta_j}^{(c_j')} \cdots e_{\beta_k}^{(c_k')} $ 
with $0 \leq c_j' \leq p^r-1$, $0 \leq c_k' \leq p^r-1$, 
and  $0 \leq c_i' \leq p^{r-1}-1$ for $j+1 \leq i \leq k-1$ and hence 
lies in $\mathcal{V}_{(r,j,k)}$. 
So we shall consider the term  
$d e_{\beta_k}^{(b_k')}  e_{\beta_s}^{(p^{r-1})}   e_{\beta_j}^{(b_j')} $ in (\ref{eqinprop2}). 
Using Proposition 
\ref{commform2} repeatedly, we see that 
$e_{\beta_k}^{(b_k')}  e_{\beta_s}^{(p^{r-1})}   e_{\beta_j}^{(b_j')} $ is 
an $\mathbb{F}_p$-linear combination of elements of the form 
$ e_{\beta_j}^{(c_j'')} \cdots e_{\beta_k}^{(c_k'')} $ with $0 \leq c_s'' \leq p^r-1$ and 
$0 \leq c_i'' \leq p^{r-1}-1$ for $i \neq s$. Therefore, we have to check that 
$d=0$ (of course in $\mathbb{F}_p$, not in $\mathbb{Z}$) or $a_s=r$ in order to show that 
$d e_{\beta_k}^{(b_k')}  e_{\beta_s}^{(p^{r-1})}   e_{\beta_j}^{(b_j')} \in \mathcal{V}_{(r,j,k)}$.    
Note that $[ e_{\beta_k},  e_{\beta_j}]=d e_{\beta_s}$ in $\mathcal{U}$ 
 by applying ${\rm Fr}^{r-1}$ to  (\ref{eqinprop1}).  

Suppose that $p=2$ and $\Phi$ is of type 
${\rm B}_l\ (l \geq 2)$, ${\rm C}_l\ (l \geq 2)$, or ${\rm F}_4$. Then 
by definition of $a_i$ we have 
\[a_i=
\left\{ \begin{array}{ll}
r & \mbox{if $\beta_i \in \mathcal{C}_0 \cup \mathcal{H}(\Phi_{\rm short}^+)$,} \\
r-1 & \mbox{otherwise}
\end{array} \right.. \]
So suppose that $\beta_s \not\in \mathcal{C}_0 \cup \mathcal{H}(\Phi_{\rm short}^+)$.  
Recall that $\beta_s=\beta_j+\beta_k$. Since $a_j=a_k=r$, both $\beta_j$ and $\beta_k$ 
must lie in $\mathcal{C}_0 \cup \mathcal{H}(\Phi_{\rm short}^+)$ and hence by Lemma 
\ref{root3} they must be short. Then by Remark \ref{root1rem} we have $d=0$. 

Suppose that $p=2$ and $\Phi$ is of type ${\rm G}_2$. Since 
\[a_i=
\left\{ \begin{array}{ll}
r & \mbox{if $\beta_i \in \{ \alpha_1, \alpha_2, \alpha_1+\alpha_2\}$,} \\
r-1 & \mbox{if $\beta_i \in \{ 2\alpha_1+\alpha_2, 3\alpha_1+\alpha_2, 3\alpha_1+2\alpha_2\}$}
\end{array} \right. \]
and $a_j=a_k=r$, 
 the possible 3-tuples 
$(\beta_j,\beta_k,\beta_s)$ are $(\alpha_2,\alpha_1, \alpha_1+\alpha_2)$ and 
$(\alpha_1+\alpha_2,\alpha_1, 2\alpha_1+\alpha_2)$. If 
$(\beta_j,\beta_k,\beta_s)= (\alpha_2,\alpha_1, \alpha_1+\alpha_2)$, we have 
$a_s=a_2=r$. If $(\beta_j,\beta_k,\beta_s)= (\alpha_1+\alpha_2,\alpha_1, 2\alpha_1+\alpha_2)$, 
we have $d=0$ by equality (\ref{commformg2}). 

Suppose that $p=3$ and $\Phi$ is of type ${\rm G}_2$. Since 
\[a_i=
\left\{ \begin{array}{ll}
r & \mbox{if $\beta_i \in \{ \alpha_1, \alpha_2, \alpha_1+\alpha_2, 2\alpha_1+\alpha_2\}$,} \\
r-1 & \mbox{if $\beta_i \in \{3\alpha_1+\alpha_2, 3\alpha_1+2\alpha_2\}$}
\end{array} \right. \]
and $a_j=a_k=r$, the possible 3-tuples 
$(\beta_j,\beta_k,\beta_s)$ are $(\alpha_2,\alpha_1, \alpha_1+\alpha_2)$,  
$(\alpha_1+\alpha_2,\alpha_1, 2\alpha_1+\alpha_2)$,  
$(\alpha_1+\alpha_2, 2\alpha_1+\alpha_2, 3\alpha_1+2\alpha_2)$, and 
$(2\alpha_1+\alpha_2, \alpha_1, 3\alpha_1+\alpha_2)$.  If 
$(\beta_j,\beta_k,\beta_s)= (\alpha_2,\alpha_1, \alpha_1+\alpha_2)$, we have 
$a_s=a_2=r$.  If $(\beta_j,\beta_k,\beta_s)= (\alpha_1+\alpha_2,\alpha_1, 2\alpha_1+\alpha_2)$, 
we have $a_s=a_4=r$.  If $(\beta_j,\beta_k,\beta_s)= 
(\alpha_1+\alpha_2,2\alpha_1+\alpha_2, 3\alpha_1+2\alpha_2)$, 
we have $d=0$ by equality (\ref{commformg5}). 
If $(\beta_j,\beta_k,\beta_s)= (2\alpha_1+\alpha_2,\alpha_1, 3\alpha_1+\alpha_2)$, 
we have $d=0$ by equality (\ref{commformg3}). 

Therefore, we obtain 
$d e_{\beta_k}^{(b_k')}  e_{\beta_s}^{(p^{r-1})}   e_{\beta_j}^{(b_j')} \in \mathcal{V}_{(r,j,k)}$ 
and hence $e_{\beta_k}^{(b_k)} e_{\beta_j}^{(b_j)} \in 
\mathcal{V}_{(r,j,k)}$. Thus the proposition is proved. \qed
\end{Proof}

Now we have completed the proof of Theorem \ref{mainthm1}. \qed \\
\ 

By symmetry, we see that a result similar  to Theorem \ref{mainthm1} also holds 
for $\mathcal{V}_r^-$ and $\mathcal{U}_r^-$. 
\ \\

\begin{Theorem}\label{mainthm1'}
Let $\Theta$ and $a_{\alpha}$ for $\alpha \in \Phi^+$  be as in Theorem \ref{mainthm1}. 
Then $\mathcal{V}_r^-$ has 
$\left\{ \left. \prod_{\alpha \in \Phi^+} e_{-\alpha}^{(n_{-\alpha})}\ \right| \  
0 \leq n_{-\alpha} \leq p^{a_{\alpha}}-1\right\}$ as an $\mathbb{F}_p$-basis, where  
an ordering of $\Phi^+$  in the product is fixed arbitrarily, and 
$\Lambda_r^- \cup \left\{ \left. e_{-\alpha}^{(p^{r-1})}\ \right| \ 
\alpha \in \Theta \right\}$ is a minimal generating set of the  
$\mathbb{F}_p$-subalgebra $\mathcal{U}_r^-$ of $\mathcal{U}$.
\end{Theorem}
\  

Finally, we deal with $\mathcal{V}_r$ and $\mathcal{U}_r$. \\
\ 

\begin{Proposition}\label{vandu}
$\mathcal{V}_r$ contains $\mathcal{U}_r^0$. 
\end{Proposition}

\begin{Proof} 
Use induction on $r$. It is easy when $r=1$, so assume that $r>1$. 
Since $\mathcal{V}_{r-1}$ contains $\mathcal{U}_{r-1}^0$ by induction, 
we only have to check that ${h_i \choose p^{r-1}} \in \mathcal{V}_r$ 
for $1 \leq i \leq l$. Note that 
\[e_i^{(p^{r-1})} f_i^{(p^{r-1})} - f_i^{(p^{r-1})} e_i^{(p^{r-1})} - 
{h_i \choose p^{r-1}} = \sum_{s=1}^{p^{r-1}-1} f_i^{(p^{r-1}-s)} 
{h_i -2p^{r-1} +2s \choose s} e_i^{(p^{r-1}-s)}\]
in $\mathcal{U}$. Since the right-hand side lies in $\mathcal{V}_{r-1}$, we see that  
${h_i \choose p^{r-1}} $ lies in $\mathcal{V}_{r}$.  \qed
\end{Proof}
\

\begin{Lemma}\label{commform3} 
For $s,t \in \mathbb{Z}_{> 0}$, 
let $x=\prod_{k=1}^s e_{i_k}^{(m_k)}$ and $y=\prod_{k=1}^t f_{j_k}^{(n_k)}$ be two elements 
in $\mathcal{U}$ with $i_k,j_k \in \{1, \dots, l\}$ and  
$m_k, n_k \in \mathbb{Z}_{>0}$. Then 
$x \mathcal{U}_r^0 =  \mathcal{U}_r^0 x$ and $y \mathcal{U}_r^0 =  \mathcal{U}_r^0 y$.  
\end{Lemma}

\begin{Proof}
The lemma follows from the fact that 
\[ {h_i+c \choose n} x = x {h_i+c + \sum_{k=1}^s \langle \alpha_{i_k}, 
\alpha_i^{\vee} \rangle m_k \choose n}\]
and
\[{h_i+c \choose n} y = y {h_i+c - \sum_{k=1}^t \langle \alpha_{j_k}, 
\alpha_i^{\vee} \rangle n_k \choose n}\]
for $0 \leq n \leq p^r-1$, $c \in \mathbb{Z}$, and $i \in \{ 1, \dots, l\}$. \qed
\end{Proof}
\ 

\begin{Lemma}\label{commform4}
 For $s,t \in \mathbb{Z}_{> 0}$, 
let $x=\prod_{k=1}^s e_{i_k}^{(m_k)}$ and $y=\prod_{k=1}^t f_{j_k}^{(n_k)}$ be two elements 
in $\mathcal{U}$ with $i_k,j_k \in \{1, \dots, l\}$, $1 \leq m_k \leq p^r-1$, and 
$1 \leq n_k \leq p^r-1$. Then the 
product $xy$ is a sum of elements of the form 
\[\left( \prod_{k=1}^t f_{j_k}^{(n_k')}\right) z \left( \prod_{k=1}^s e_{i_k}^{(m_k')}\right)\]
with $1 \leq m_k' \leq p^r-1$, $1 \leq n_k' \leq p^r-1$, and $z \in \mathcal{U}_r^0$. 
\end{Lemma}

\begin{Proof}
Use induction on $s,t$ together with Lemma \ref{commform3}. 
\end{Proof}
\ 

Now we obtain the second main result. \\

\begin{Theorem}\label{mainthm2} 
Let $\Theta$ and $a_{\alpha}$ for $\alpha \in \Phi^+$  be as in Theorem \ref{mainthm1}. 
The following hold. \\

\noindent {\rm (i)} The multiplication map 
$\mu : \mathcal{V}_r^- \otimes_{\mathbb{F}_p} \mathcal{U}_r^0 \otimes_{\mathbb{F}_p}  
 \mathcal{V}_r^+ \rightarrow \mathcal{V}_r$ is an $\mathbb{F}_p$-linear isomorphism. 
In particular, $\mathcal{V}_r$ has 
\[\left\{ \left. \prod_{\alpha' \in \Phi^-} e_{\alpha'}^{(n_{\alpha'}')}
\prod_{i=1}^l {h_i \choose m_i}   
\prod_{\alpha \in \Phi^+} e_{\alpha}^{(n_{\alpha})} \ \right| \  
\begin{array}{l}
{0 \leq n_{\alpha} \leq p^{a_{\alpha}}-1, 0 \leq n_{\alpha'}' \leq p^{a_{-\alpha'}}-1,} \\
{0 \leq m_i \leq p^r-1} 
\end{array}
\right\}\] 
as an $\mathbb{F}_p$-basis, where  
 orderings of $\Phi^+$ and $\Phi^-$ in  the products are fixed arbitrarily. \\

\noindent {\rm (ii)} 
$\Lambda_r \cup \left\{ \left. e_{\alpha}^{(p^{r-1})}, e_{-\alpha}^{(p^{r-1})}\  \right| \ 
\alpha \in \Theta \right\}$ generates the $\mathbb{F}_p$-algebra $\mathcal{U}_r$. 
\end{Theorem}

\begin{Proof}
(ii) follows from Theorems \ref{mainthm1} and \ref{mainthm1'}. 
Let us  prove  (i). Proposition \ref{vandu} implies that the map $\mu$ is well-defined. 
Since ${\rm Im} \mu$ contains $\Lambda_r$, we only have to  
show that the image is closed under multiplication. 
But it follows easily from Lemmas \ref{commform3} and \ref{commform4}. \qed
\end{Proof}
\

\noindent {\bf Acknowledgements.} \ \ 
The author would like to thank the referee for carefully
reading the manuscript and giving some helpful comments. 
This work was supported by JSPS KAKENHI Grant Number JP18K03203.



\end{document}